\pdfoutput=1
\documentclass[journal]{IEEEtran}
\ifCLASSINFOpdf
\else
\fi
\usepackage{graphicx}
\usepackage{caption}
\usepackage{subcaption}
\graphicspath{{./figs/}}
\usepackage{color} 

\usepackage{import}
\usepackage{epsfig}

\usepackage{paralist}

\usepackage{amssymb}
\usepackage{amsmath}
\usepackage{bm}

\usepackage[mathscr]{euscript}
\usepackage{comment}

\usepackage[hyphens]{url}
\usepackage[capitalise,english]{cleveref} 

\usepackage{framed} 
\usepackage{multicol} 

\usepackage{multirow} 

\usepackage[noprefix]{nomencl} 
\makenomenclature
\renewcommand*\nompreamble{\begin{multicols}{2}}
\renewcommand*\nompostamble{\end{multicols}}

\usepackage{numprint}
\npthousandsep{,}\npthousandthpartsep{}\npdecimalsign{.}

\usepackage{lineno}

\usepackage{xcolor}
\colorlet{lightred}{red!70!white}

\usepackage[noadjust]{cite}

\newcommand{\eg}{e.g.} 
\newcommand{\ie}{i.e.} 

\makeatletter
\newcommand*{\rom}[1]{\expandafter\@slowromancap\romannumeral #1@}
\makeatother

\begin{document}
%

\title{Real-time Error Control for Surgical Simulation}

\author{\IEEEauthorblockN{Huu Phuoc Bui\IEEEauthorrefmark{1}\IEEEauthorrefmark{2}, Satyendra Tomar\IEEEauthorrefmark{2}, Hadrien Courtecuisse\IEEEauthorrefmark{1}, St\'{e}phane Cotin\IEEEauthorrefmark{3} and St\'{e}phane Bordas\IEEEauthorrefmark{2}\IEEEauthorrefmark{4}}

\IEEEauthorblockA{\IEEEauthorrefmark{1}University of Strasbourg, CNRS, ICube, F-67000 Strasbourg Fance\\
\IEEEauthorrefmark{2}University of Luxembourg, Research Unit of Engineering Science, L-1359 Luxembourg, Luxembourg\\
\IEEEauthorrefmark{3}Inria Nancy Grand Est, 54603 Villers-les-Nancy, France\\
\IEEEauthorrefmark{4}Cardiff University, School of Engineering, Queen’s Buildings, The Parade, Cardiff CF24 3AA, Wales, UK
}
}

\maketitle


\begin{abstract}
\emph{Objective}: To present the first real-time a posteriori error-driven adaptive finite element approach for real-time simulation and to demonstrate the method on a needle insertion problem. \emph{Methods}: We use corotational elasticity and a frictional needle/tissue interaction model. The problem is solved using finite elements within SOFA\footnote{https://www.sofa-framework.org/}. The refinement strategy relies upon a hexahedron-based finite element method, combined with a posteriori error estimation driven local $h$-refinement, for simulating soft tissue deformation. \emph{Results}: We control the local and global error level in the mechanical fields (\eg{} displacement or stresses) during the simulation.  We show the convergence of the algorithm on academic examples, and demonstrate its practical usability on a percutaneous procedure involving needle insertion in a liver. For the latter case, we compare the force displacement curves obtained from the proposed adaptive algorithm with that obtained from a uniform refinement approach. \emph{Conclusions}: Error control guarantees that a tolerable error level is not exceeded during the simulations. Local mesh refinement accelerates simulations. \emph{Significance}: Our work provides a first step to discriminate between discretization error and modeling error by providing a robust quantification of discretization error during simulations.
      
\end{abstract}


%
\IEEEpeerreviewmaketitle

\section{Introduction}
\label{sec:introduction}

\subsection{Motivation}
\IEEEPARstart{R}{}eal-time simulations are becoming increasingly common for various applications, from geometric design \cite{nealen2006physically,wang2015linear} to medical simulation \cite{Courtecuisse2014}. Our focus is on real-time simulation of the interaction of a surgeon or interventional radiologist with deformable organs. Such simulations are useful to both, help surgeons train, rehearse complex operations or\slash and to guide them during the intervention. In time, reliable simulations could also be central to robotic surgery.

A number of factors are concurrently involved in defining the ``accuracy’’ of surgical simulators: mainly the modeling error and the discretization error. Most work in the area has been looking at the above sources of error as a compounded, lumped, overall error. Little or no work has been done to discriminate between modeling error (e.g. needle-tissue interaction, choice of constitutive models) and discretization error (use of approximation methods like FEM). However, it is impossible to validate the complete surgical simulation approach and, more importantly, to understand the sources of error without evaluating both the discretization error and the modeling error.

The first ingredient in any mechanical simulation is the ability to simulate the deformation of the solid of interest. This deformable solid mechanics problem is usually solved by finite element method (FEM) \cite{zienkiewicz1977finite} or meshless\slash meshfree methods \cite{nguyen2008meshless}, which are used to discretize the equilibrium equations. It is usually uneconomical or prohibitively expensive to use a fixed mesh for such simulations. Indeed, coarse meshes are sufficient to reproduce ``smooth'' behavior, whereas ``non-smooth'' behavior such as discontinuities engendered by cuts or material interfaces, singularities, boundary layers or stress concentrations require a finer mesh. Adequate approaches are thus needed to refine the discretization in these areas. 

Yet, existing numerical methods used in surgical simulation use either a fixed discretization (finite element mesh, meshfree point cloud, reduced order method), or adapt the mesh using heuristics \cite{wu2001adaptive,muller2002stable}. To our knowledge, no approach is currently able to adapt the finite element mesh based on rational a posteriori error estimates \cite{ainsworth2011posteriori}.

Our objective is thus to devise a robust and fast approach to local remeshing for surgical simulations. To ensure that the approach can be used in clinical practice, the method should be robust enough to deal, as realistically as possible, with the interaction of surgical tools with the organ, and fast enough for real-time simulations. The approach should also lead to an improved convergence so that an ``economical'' mesh is obtained at each time step. The final goal is to achieve optimal convergence and the most economical mesh, which will be studied in our future work.

In this paper, we propose and benchmark a local mesh refinement and coarsening approach which is based on the estimation of the discretization error incurred by FEM in the solution of a corotational model representing soft tissues.  The general ideas presented here can be used directly in geometric design based on deformable models. Our proposed approach has similitudes with the octree approaches of \cite{seiler2011robust}.

\subsection{Error in numerical simulations}

It is useful to first review the various sources of error in numerical simulations. The first error source arises when a mathematical model is formulated for a given physical problem: this is known as the modeling error. The second error arises upon discretization of this mathematical model, for example using FEM or meshfree methods. Finally, numerical error is incurred because of the finite precision of computers and round off errors. In this paper, we focus on the second source of error, namely discretization error. We therefore assume that the model we use is descriptive of reality, \ie{} we are solving the right problem, and we ask ourselves the question whether ``we are solving the problem right'', in other words, correctly. 

The main difficulty in answering this question comes from the fact that an exact solution, to which the numerical solution could be compared, is generally not available. Different approaches exist to address this problem, which are reviewed in the literature, see \eg{}, \cite{zienkiewicz1977finite, ainsworth2011posteriori, Verfuerth-13-Apost}. Simple methods available in practice to indicate the error distribution can be categorized into two classes: recovery-based and residual-based.

The first class of indicators assumes that the exact solution (of stresses) is smooth enough (at least locally). They rely on the construction of an ``improved'' numerical solution from the raw numerical solution, to which the raw numerical solution can be compared. Where these two solutions are significantly different (above certain threshold), the error level is high and the mesh should be refined, and where these two solutions are close together, the mesh can be kept unchanged or coarsened. This idea was proposed by Zienkiewicz and Zhu in \cite{zienkiewicz1992superconvergent}, and its asymptotic convergence to the exact error is studied in \cite{CarstensenB-02-Apost, BartelsC-02-Apost}.

The second class of indicators relies on the computation of the residual of the governing equations within each computational cell, typically each finite element. These ``residual-based'' error indicators lead to mesh refinement where the solution leads to large residuals, and keep the mesh constant or coarsen it where the element residuals are relatively small, compared to a given tolerance. These estimates were first proposed by Babuska and Rheinboldt in \cite{BabuskaR-78-Apost}.

Based on these local error estimators, mesh refinement methodologies can be devised to derive mesh adaptation, see \eg{}~\cite{Verfuerth-13-Apost} for a recent comprehensive presentation. This requires two key ingredients: \emph{a marking strategy} that decides which elements should be refined, and \emph{a refinement rule} that defines how the elements are subdivided. For element marking, we use the \emph{maximum strategy}, see \cref{sec:Error_estimate_and_adaptive_refinement} for details. Other strategies, such as bulk\slash equilibration strategy or percentage strategy, see \eg{}~\cite{Verfuerth-13-Apost}, can also be used.

\subsection{Simulation of percutaneous operations}

Needle-based percutaneous procedures are an important part of modern clinical interventions such as biopsy, brachytherapy, cryotherapy or regional anesthesia. The success of these procedures depends on good training and careful planning to optimize the path to the target, while avoiding critical structures \cite{Hamze2016}. In some instances the procedure can also be assisted by robotic devices. Unfortunately, natural tissue motion (due to breathing, for instance), and deformation (due to needle insertion) generally lead to incorrect or inefficient planning \cite{Hamze2016}. To address these issues, one must rely on an accurate simulation of needle insertion. For most problems, computational speed is also very important, since the simulation is at the core of an optimization algorithm (for the needle path) or a robotic control loop.

The main works on needle insertion (see the survey by \cite{Abolhassani07}) propose to model the interaction between the needle and soft tissues using FEM. In the various methods proposed in the literature, three main research directions have been followed: soft tissue model, flexible needle model and needle-tissue interactions. The needle model is usually not an issue, both in terms of modeling choice and computational cost. For instance, in \cite{Duriez2009}, authors report computation times of a few milliseconds for a FEM needle model composed of 50 serially-linked Timoshenko beam elements. Soft tissue models are usually based on FEM, and rely on linear or non-linear constitutive laws \cite{Abolhassani07}.

A large body of work covers the modeling and simulation of soft tissue deformation, even under real-time computation constraints. But overall, the interaction model between the needle and tissue remains a major  challenge. It combines different physical phenomena, such as puncturing, cutting, sliding with friction and Poynting's effect. To capture the essential characteristics of these interactions, existing methods usually rely on experimental force data and remeshing techniques in order to align nodes of the FEM mesh with the needle path. In \cite{Duriez2009}, a constraint-based approach, avoiding remeshing, was used to simulate needle-tissue interactions. However, the simulations did not account for realistic anatomical details. In addition, Misra {\em et al.} \cite{Misra2010} showed that needle steering, which occurs when using asymmetric needle tips, can be modeled using microscopic observations of needle-tissue interactions. 
 
Unfortunately, if no assumption can be made about the region of the domain where the needle will be inserted, simulations involving very detailed meshes become very slow, which is a real issue in the context presented above.
Error-controlled real-time simulation of needle insertion is thus an unsolved problem whose solution requires tackling a number of difficulties:
\begin{itemize}
\item developing needle-tissue models. A review of cutting simulation is provided in \cite{wu2014physically}
\item using these models within discrete approaches like FEM, mesh-free methods, or others
\item accelerating the simulation (advanced hardware, model order reduction)
\item validating the needle-tissue interaction model combined with discrete solution (are we solving the right problem?)
\item verifying the discrete solution, i.e. controlling the discretization error associated with the discrete model (are we solving the problem right?)
\end{itemize}

In this paper, we propose to focus on the last point above, with the aim to model needle-tissue interactions using an adaptive meshing strategy driven by simple a posteriori error estimation techniques. Similarly to \cite{Duriez2009}, we do not require the mesh to conform to the needle path. Mesh subdivision is only introduced as a means to improve the accuracy of the needle-tissue interactions. Our mesh refinement method is guided by the stress\slash energy error estimate resulting from the needle-tissue interaction and imposed boundary conditions: elements of the mesh are subdivided when a numerical error threshold is reached. The subdivision process is completely reversible, i.e. refined elements are set back to their initial topology when refinement is no longer needed. Our refinement approach does not rely on a usual octree structure (see also \cite{dick2011hexahedral}), thus allowing a variety of subdivision schemes that are well suited for needle insertions, as detailed in \cref{sec:Error_estimate_and_adaptive_refinement}. Using this approach, interactive computation times can be achieved while detailed tissue motion near the needle shaft or tip can be computed. This opens new possibilities for fast simulations of flexible needle insertion in soft tissues. We illustrate the convergence study of the adaptive refinement scheme and some of the possible scenarios in \cref{sec:results}.


\section{Model and discretization}
\label{sec:Model_and_discretization}

In this section, we describe the model and the discretization approach, which are used for needle and soft tissue interaction.
\subsection{Problem statement}
\begin{figure}[!htbp]
\centering
\includegraphics[width=0.55\columnwidth]{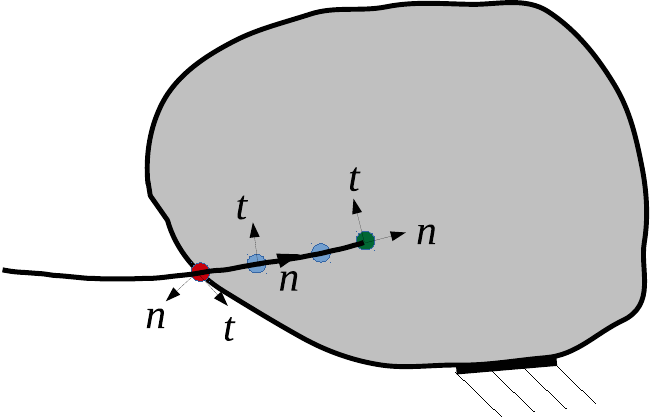}
\caption{Three types of constraints between the needle and soft tissue: surface puncture (in red), needle tip constraint (in green) and needle shaft constraints (in blue). A local coordinate system $n$-$t$ is defined at each constraint point.}
\label{fig:needle_tissue_schema}
\end{figure}

During the needle insertion, three types of constraints are defined, see \cref{fig:needle_tissue_schema}.
Coulomb's friction law is used to describe frictional contacts within these three types of constraints. First, a puncture constraint is defined between the needle tip and the tissue surface. This constraint satisfies the Kuhn-Tucker condition in the direction $n$ (normal to the tissue surface)
\begin{equation}
\delta_n \geq 0, \qquad \lambda_n \geq 0, \qquad \delta_n \cdot \lambda_n = 0,
\end{equation}
where $\delta_n$ denotes the distance between the needle tip and the tissue surface in the direction $n$, and $\lambda_n$ denotes the contact force in that direction. Let $\lambda_{n0}$ represent the puncture strength of the tissue. The Kuhn-Tucker condition expresses that the contact force only exists when the needle tip is in contact with the tissue surface. When the contact force is higher than a threshold (puncture strength of the tissue)
\begin{equation}
\lambda_n > \lambda_{n0},
\end{equation}
the needle can penetrate into the tissue. In the tangent direction $t$, Coulomb's friction law is considered in order to take into account the stick\slash slip between the needle tip and tissue surface
\begin{equation}
\lambda_t < \mu \lambda_n  \quad \textrm{(stick)}; \quad \lambda_t = \mu \lambda_n  \quad \textrm{(slip)},
\end{equation}
where $\mu$ denotes the friction parameter.

Second, a needle tip constraint is defined at the tip of the needle as soon as it penetrates into the tissue. Depending on the relationship between the contact forces in the normal direction $n$ (along the needle shaft) and in the tangent direction $t$ (see \cref{fig:needle_tissue_schema}),
the needle tip can cut and go through the tissue or not
\begin{equation}
  \lambda_n < \mu \lambda_t + \lambda_{n0} \; \textrm{(stick)}; \quad \lambda_n \geq  \mu \lambda_t + \lambda_{n0}  \; \textrm{(cut and slip).}
 \end{equation}

Finally, needle shaft constraints are defined along the needle shaft so that the needle shaft is enforced to follow the insertion trajectory created by the advancing needle tip. Again, the Coulomb's friction law is applied to these constraints to represent the stick and sliding contact between the tissue and the needle shaft
\begin{equation}
\lambda_n  < \mu \lambda_t \; \textrm{(stick)}; \quad \lambda_n  = \mu \lambda_t \; \textrm{(sliding).}
\end{equation}

\subsection{Strong form}

\begin{figure}[!htbp]
 \centering
\includegraphics[width=0.9\columnwidth]{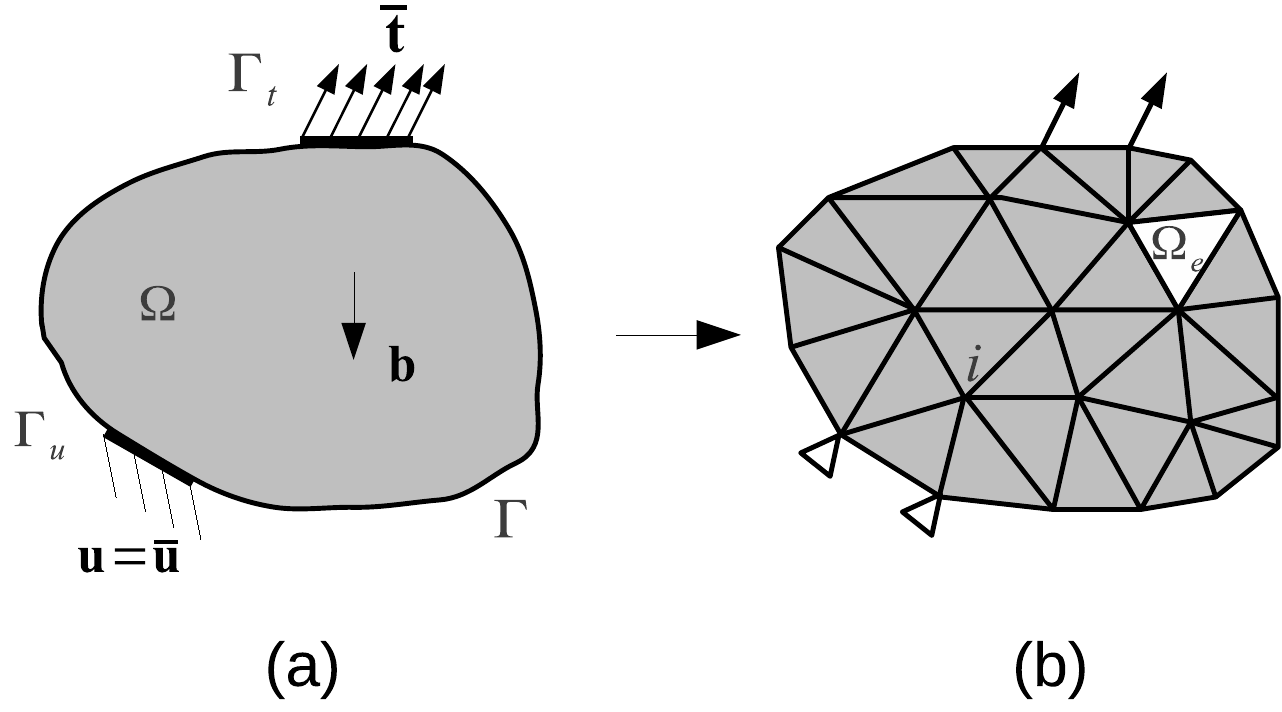}
\caption{A body $\Omega$ subjected to a traction $\mathbf{\bar{t}}$ on its boundary part $\Gamma_t$, a body force $\mathbf{b}$, and an imposed displacement $\mathbf{\bar{u}}$ on boundary part $\Gamma_u$ (a); Simplified illustration of FEM discretization (b).}
\label{fig:fem_discretization}
\end{figure}

We model both, the tissue and the needle, as dynamic deformable objects. Thus, they can be regarded as dynamic elastic solids, and the governing equations of the model are formulated as
\begin{subequations}
\begin{align}
 \mathrm{div}\, \bm{\sigma} + \mathbf{b} + \bm{\lambda} & =  \rho \ddot{\mathbf{u}} \qquad \text{in } \Omega \label{eq:equilibrium}\\
 \bm{\epsilon} & =  \frac{1}{2}\left(\mathrm{grad}\,\mathbf{u} + (\mathrm{grad}\,\mathbf{u})^T \right)  \label{eq:kinematic} \\
 \bm{\sigma} & =  f(\bm{\epsilon},\bm{\nu}) \label{eq:constitutive}\\
 \bm{\sigma} \cdot \mathbf{n} & = \mathbf{\bar{t}} \; \text{on } \Gamma_t \label{eq:NeumannBC}; \\
 \mathbf{u}  & = \mathbf{\bar{u}} \; \text{on } \Gamma_u, \label{eq:DirichletBC}
\end{align}
\end{subequations}
where $\bm{\sigma}$ is the Cauchy stress tensor, $\mathbf{b}$ is the body force vector, $\rho$ is the mass density, $\bm{\epsilon}$ is the strain tensor, $\bm{\nu} = (\nu_1, \nu_2, \dots, \nu_n)$ is the internal variables, (~$\dot{}$~) denotes the partial derivative with respect to time, $\mathbf{n}$ denotes the outward unit normal vector on $\Gamma_t$, and $\bm{\lambda}$ denotes the contact force between the needle and the tissue. The object domain and boundary conditions are shown in~\cref{fig:fem_discretization}a.

\subsection{Spatial and temporal discretization} 

%
\subsubsection{Space discretization}


The basic idea of FEM is to discretize the domain $\Omega$ into finite elements $\Omega_e$, $e=1,2,\dots,N_e$, by $N_n$ nodes, as depicted in Figure~\ref{fig:fem_discretization}b. Based on the discretization concept, see \eg{}, \cite{Liu201443,zienkiewicz2000finite}, we obtain the discrete problem of the element $e$
\begin{equation}
\label{eq:elmentEquation}
 \mathbf{M}_e \ddot{ \mathbf{u} }_e + \mathbf{C}_e \dot{ \mathbf{u} }_e + \mathbf{f}_e (\bm{\sigma}) = \mathbf{f}_e^{ext},
\end{equation}
where $\mathbf{M}_e$ is the element mass matrix, $\mathbf{K}_e$ is the element stiffness matrix, $\mathbf{C}_e$ is the damping matrix, $\mathbf{f}_e^{ext}$ is the external force applied to the element $e$, while the internal force $\mathbf{f}_e (\bm{\sigma})$ reads
\begin{equation}
\begin{split}
\label{eq:f=Ku_linear}
 \mathbf{f}_e (\bm{\sigma}) = \int_{\Omega_e} \mathbf{B}_e^T \bm{\sigma}  \, \mathrm{d}\Omega & = \int_{\Omega_e} \mathbf{B}_e^T \mathbf{E} \mathbf{B}_e \mathbf{u}_e   \, \mathrm{d}\Omega \\
 & = \mathbf{K}_e \cdot \mathbf{u}_e = \mathbf{K}_e \cdot (\mathbf{x}_e - \mathbf{x_0}_e)
 \end{split}
\end{equation}
where $\mathbf{B}_e$ is the strain-displacement matrix, $\mathbf{E}$ is the fourth-order stiffness tensor, $\mathbf{x}_e$ and $\mathbf{x_0}_e$ denote the current and initial position of the element $e$, respectively.
However, using~\cref{eq:f=Ku_linear} results in inaccuracy for large rotations problems, which is observed by artificially inflated deformation of the elements. To overcome this, Felippa \textit{et al.}~\cite{Felippa2005} decomposed the deformation gradient of the element into the rigid and deformation parts, and the element nodal internal force becomes
\begin{equation}
 \mathbf{f}_e = \mathbf{R}_e \mathbf{K}_e ( \mathbf{R}^T_e \mathbf{x}_e - \mathbf{x_0}_e),
\end{equation}
where $\mathbf{R}$ stands for the element rotation matrix of the element local frame with respect to its initial orientation, being updated at each time step. Using this corotational formulation results in no visual artifacts.

The global mass, stiffness and damping matrices of the system can then be assembled from the element ones, and \cref{eq:elmentEquation} can be rewritten to a global system equation as
\begin{equation}
\label{eq:Ma=f}
 \mathbf{M} \mathbf{a} = \mathbf{f}(\mathbf{x},\mathbf{v}),
\end{equation}
where $\mathbf{a} = \ddot{\mathbf{u}}$, $\mathbf{x}$, $\mathbf{v} = \mathbf{\dot{\mathbf{u}}}$ are the acceleration, position and velocity vectors, respectively, and $\mathbf{f}(\mathbf{x},\mathbf{v}) = \mathbf{f}^{ext}-\mathbf{K}\mathbf{u} - \mathbf{C}\mathbf{v}$ represents the net force (the difference of the external and internal forces) applied to the object.

In our simulations, a diagonally lumped mass matrix is employed, and the stiffness matrix $\mathbf{K}$ is computed based on the co-rotational FE formulation described above, which allows large rotations for both, needle as well as tissue. For higher accuracy of the computed strain field, the soft tissue domain is discretized using hexahedral elements. To avoid the complex issue of generating an exact hexahedral mesh of the domain, we use a mesh that does not conform to the boundary of the domain, as in Immersed Boundary Method~\cite{Pinelli2010}.
The needle, on the other hand, is modeled using serially-linked beam elements, as in~\cite{Duriez2009}. In this case, each node of the needle has $6$ degrees of freedom ($3$ translations and $3$ rotations), while the tissue model only uses $3$ translational degrees of freedom per node.


Since the FEM formulation is based on the discretization of the physical domain, it naturally introduces the discretization error in the result. To control this error source, in \cref{sec:Error_estimate_and_adaptive_refinement} we present an adaptive refinement scheme.

\subsubsection{Time discretization}
For temporal discretization, we use an implicit backward Euler scheme~\cite{Baraff1998}, which is described as follows
\begin{equation}
\label{eq:EulerBackward}
\dot{\mathbf{u}}_{t+\tau} = \dot{\mathbf{u}}_t + \tau\ddot{\mathbf{u}}_{t+\tau}; \quad \mathbf{u}_{t+\tau} = \mathbf{u}_{t} + \tau \dot{\mathbf{u}}_{t+\tau},
\end{equation}
where $\tau$ denotes the time step. Inserting \cref{eq:EulerBackward} into \cref{eq:Ma=f} gives the final discrete system
\begin{equation}
\label{eq:Ax=b}
\underbrace{(\mathbf{M}-\tau \mathbf{C}-\tau^2 \mathbf{K})}_{\mathbf{A}} d\mathbf{v} = \underbrace{ \tau\mathbf{f}(\mathbf{x}^t,\mathbf{v}^t) + \tau^2 \mathbf{K}\mathbf{v}^t}_{\mathbf{b}}
\end{equation}
where $d\mathbf{v} = \mathbf{v}_{t+\tau} - \mathbf{v}_t$.

After solving~\eqref{eq:Ax=b} for $d\mathbf{v}$, the position and velocity are updated for needle and tissue as
\begin{equation}
\label{eq:updateVX}
\mathbf{v}_{t+\tau} = d\mathbf{v} + \mathbf{v}_t; \quad \mathbf{x}_{t+\tau} = \mathbf{x}_{t} + \tau \mathbf{v}_{t+\tau}.
\end{equation}

\subsection{Constraint enforcement for needle-tissue interaction} 
To take into account the needle-tissue interaction when they are in contact, a constrained dynamic system is solved for the needle and the tissue. \cref{eq:Ax=b} then becomes
\begin{equation}
\mathbf{A} d\mathbf{v} = \mathbf{b} + \mathbf{J} \bm{\lambda},
\end{equation}
where $\bm{\lambda}$ denotes Lagrange multipliers representing the interaction forces between the needle and the tissue, and $\mathbf{J}$ provides the direction of the constraints. Different types of constraints between needle and tissue are used, and solving their interaction is detailed in \cref{sec:needle_tissue_interaction_algo}.

\textbf{Remark} Combining more advanced and clinically relevant needle-tissue interaction is straightforward in our approach.


\section{Error estimate and adaptive refinement}
\label{sec:Error_estimate_and_adaptive_refinement}

To achieve faster and more accurate FEM simulations, different adaptive techniques have been proposed in the literature. Octree-based approaches~\cite{Kwak2003531} are the most common, but the refinement procedure is limited to cubic elements which are recursively subdivided into eight finer elements. To overcome this limitation, more generic remeshing techniques~\cite{Koschier2014,Burkhart2010,Paulus2015} have been proposed. However, they are complex to implement, and may lead to ill-shaped elements. Our template-based refinement algorithm is designed to be independent of the type of element (tetrahedra, hexahedra, others), and produces a high quality mesh (thanks to the well-shaped elements of the predefined template).

Starting with an initial, relatively coarse mesh (as required to achieve real-time simulation), a criterion based on a posteriori error estimate is evaluated to drive the local refinement. The elements where the stress increases, \ie{} $ {\mathrm d \bm{\sigma}}/{\mathrm d t}>0 $, are considered for refinement, and the elements where the stress decreases, \ie{} ${\mathrm d \bm{\sigma}}/{\mathrm d t}<0$, are taken to a lower refinement (coarsening) level.
We define the approximate error of an element $\Omega_{e}$ as
\begin{equation}
\label{eq:energyError}
\eta_{e} = \sqrt{\int_{\Omega^e} (\bm{\epsilon}^h - \bm{\epsilon}^s)^T (\bm{\sigma}^h - \bm{\sigma}^s)  \mathrm d \Omega },
\end{equation}
%
which is the energy norm of the distance between the FEM solution (denoted by $h$) and an improved solution (denoted by $s$) obtained by the Zienkiewicz-Zhu smoothing procedure~\cite{zienkiewicz1992superconvergent}.
Among the elements with increasing stress, only those elements, where the error exceeds the predefined threshold, are subdivided (refined).
Similarly, among the elements with decreasing stress, only those elements, where the error is smaller than the above threshold, are coarsened.
Notice that we are not limited by the regularity of the mesh and can start from any (reasonable) heterogeneous mesh as a starting point prior to refinement.

\subsection{Zienkiewicz-Zhu error estimate}
\label{sec:Zienkiewicz-Zhu}

\begin{figure}[!htbp]
\centering
\includegraphics[width=0.5\columnwidth]{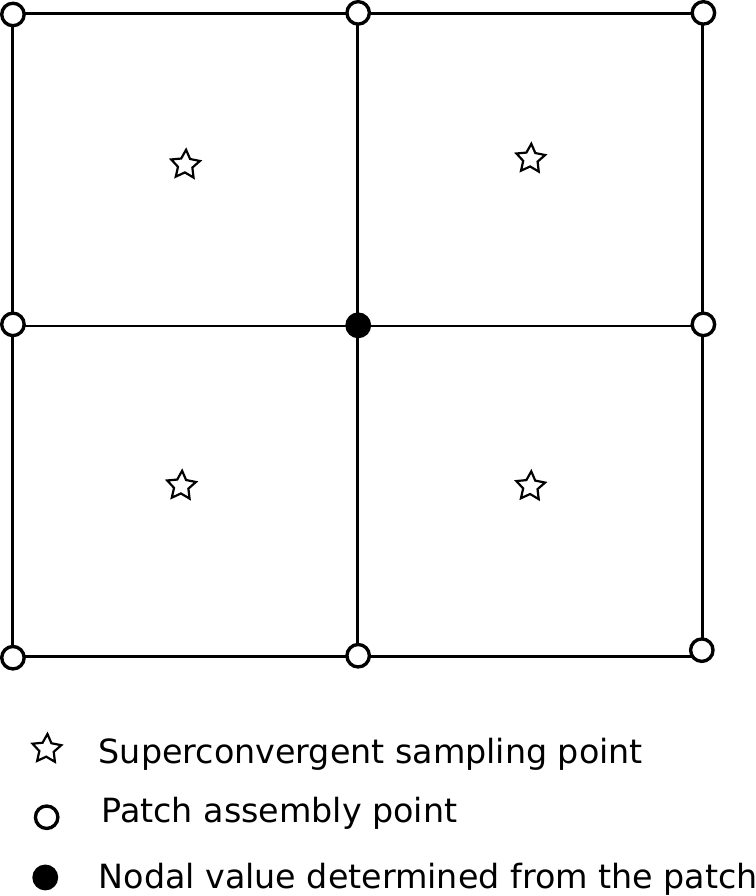}
\caption{The smoothed gradient is obtained from an element patch.}
\label{fig:SPR}
\end{figure}

Using the superconvergent patch recovery (SPR) procedure~\cite{zienkiewicz1992superconvergent}, the smoothed stress field $\bm{\sigma}^s$ is recovered from the stresses computed at the element center. The idea of this technique is based on the fact that the stress and strain at the superconvergent points (at element center in the case of linear hexahedral elements) are accurate with higher order than at the element nodes, and these values are employed to recover the nodal stress and strain within the least squares sense. A 2D representation of a patch of 8 hexahedral elements is shown in~\cref{fig:SPR}. For each component $\sigma^h_j$ of the FEM solution $\bm{\sigma}^h$, the nodal recovered stresses are computed by defining a polynomial interpolation within the element patch as
\begin{equation}
\label{eq:SPR_interpolation}
\resizebox{0.89\columnwidth}{!}{$\tilde{\sigma}_j^s = \mathbf{P} \mathbf{a}_j = [1\, x\, y\, z\, xy\, yz\, zx\, xyz] [a_1\, a_2\, a_3\, a_4\, a_5\, a_6\, a_7\, a_8]^T_j$.}
\end{equation}
Let $\mathbf{P}_k \equiv \mathbf{P}(x_k,y_k,z_k)$. To determine the unknowns $\mathbf{a}_j$ we minimize, for 8 sampling points $k$ of an element patch,
\begin{equation}
\Pi = \sum_{k=1}^8 \left[  \sigma^h_j(x_k,y_k,z_k) - \mathbf{P}_k \mathbf{a}_{j_k} \right ].
\end{equation}
This minimization results in finding $\mathbf{a}_j$ by
\begin{equation}
 \mathbf{a}_j = \mathbf{A}^{-1} \mathbf{b}
\end{equation}
where
\begin{equation}
 \mathbf{A} = \sum_{k=1}^8 \mathbf{P}_k^T \mathbf{P}_k \; \mathrm{and} \; \mathbf{b} = \sum_{k=1}^8 \mathbf{P}_k^T \sigma_j^h (x_k,y_k,z_k)
\end{equation}
Once $\mathbf{a}_j$ is available, the nodal recovered stress values are obtained by simply employing~\cref{eq:SPR_interpolation} with $\mathbf{P}$ evaluated at the corresponding node.
\footnote{Similar to the displacements, the recovered stresses $\bm{\sigma}^s$ can also be obtained using element shape functions $\bm{\sigma}^s = \mathbf{N} \bm{\tilde{\sigma}^s }$.}
%

\subsection{Element marking strategy}
\label{sec:markingStrategy}
After obtaining error distribution across all the elements, we employ the maximum strategy to select those elements which must be refined for the next level mesh. In this strategy, only those elements, where error (see \cref{eq:energyError}) is higher than certain threshold, are refined. Let $\eta_{M} = \max_{e} \eta_{e}$, where $\eta_{e}$ is defined in \cref{eq:energyError}. We mark an element for refinement if
\begin{equation}
\eta_{e} \ge \theta \eta_{M} \; \textrm{with} ~ 0< \theta <1.
\label{eq:max_rule}
\end{equation}
%

Other marking strategies, such as bulk\slash equilibration strategy, or percentage strategy, see \eg{}~\cite{Verfuerth-13-Apost}, can also be used. However, the maximum strategy described above is the cheapest among all, and hence, it is preferred for our use. In the maximum strategy, a large value of $\theta$ leads to small number of elements marked for refinement, and small value of $\theta$ leads to large number of elements marked for refinement. In our studies presented in \cref{sec:results}, we set $\theta = 0.3$.

\subsection{Template-based adaptive h-refinement}

\begin{figure}[!htbp]
\centering
\def\svgwidth{0.99\columnwidth}
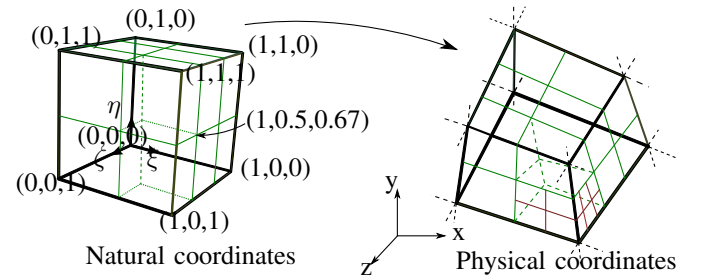
\caption{Adaptive subdivision process: each element to be subdivided is topologically transformed in its reference shape, using a template expressed in natural coordinates. Cartesian coordinates in the mesh are computed using the element shape functions. The process can be applied recursively, and is completely reversible.}
\label{fig:natural_to_cartesian_hexa}
\end{figure}

Once the refinement criterion is satisfied within an element, the element is replaced by several elements according to a predefined template. The template is simply a set of nodes and an associated topology, defined using an isoparametric formulation. The template nodes are added by using their natural coordinates.
The position $\mathbf{x}_j$, in Cartesian coordinates, of the new node $j$ is defined as $ \mathbf{x}_j = \mathbf{x}_i \xi_i^j $, where the Einstein summation convention is applied on the nodes $i$ of the removed element ($i={1,2,\dots 8}$ for hexahedral elements). The shape function $\xi_i^j$ is computed from the barycentric coordinates of  template node $j$ with respect to the node $i$. The procedure is summarized below:
{\ttfamily \small
\begin{enumerate}
\item Remove the element to be refined
\item Add template nodes and then template elements using the element shape functions  
\item Update the topology of the global mesh
\item Compute stiffness matrix of new elements
\item If needed, update the mass and damping matrices
\end{enumerate}
}
It is worth mentioning that if, after refinement, a new element fulfills the refinement criterion, it can be refined again, using the same predefined template. This results in a multi-resolution mesh (see \cref{fig:natural_to_cartesian_hexa}). Conversely, if the coarsening criterion is satisfied in already refined elements, the coarsening procedure is applied by simply removing respective fine elements, and updating the associated matrices.





\subsection{T-junction handling}
Since elements are refined by using templates, regardless of their neighboring elements, some T-junctions (incompatible nodes or hanging nodes) are generated. To avoid discontinuities at the T-junctions during the simulation, these nodes need to be handled in a special way, in which T-junction nodes are considered as slave of other independent (master) degrees of freedom (DOFs). One of the possible options is to use Lagrange multipliers, but this approach increases the total number of DOFs (as to solve the unknown Lagrange multipliers in addition), and usually leads to ill-conditioned systems. In our approach, we follow the method proposed by~\cite{Sifakis2007}, which considers only the reduced system (without T-junctions) when solving for the new positions.
Let $\mathbf{T}$ denote the transformation matrix from the reduced system to the full one (with T-junctions). The matrix $\mathbf{T}$ contains the barycentric coordinates of the T-junctions (slaves) with respect to their masters, and contains $1$ for all other normal DOFs.
The reduced system matrix $A_r$ is then computed from the full system matrix $A_f$ as
\begin{equation}
\mathbf{A}_r = \mathbf{T}^T \mathbf{A}_f \mathbf{T}.
\end{equation}
The nodal forces in the reduced space is computed from the full space as $\mathbf{f}_r = \mathbf{T}^T \mathbf{f}_f$. The reduced system $\mathbf{A}_r d\mathbf{v}_r = \mathbf{f}_r$ is solved to find $d\mathbf{v}_r$, the difference of velocity between current and previous time step.  Once $d\mathbf{v}_r$ is available, the difference of velocity in the full space are easily deduced as $d\mathbf{v}_f = \mathbf{T} d\mathbf{v}_r$. The latter is employed to update the new position and velocity of the object as in \cref{eq:updateVX}.

\begin{figure}[!htbp]
\centering
\def\svgwidth{0.6\columnwidth}
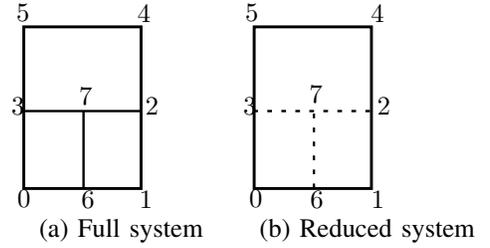
\caption{Illustration of T-junction handling method on a schematic example.}
\label{fig:t_junction_schema}
\end{figure}

A heuristic example is shown in~\cref{fig:t_junction_schema} to explicitly illustrate the method for T-junction handling, especially how the full and reduced systems are defined, resulting in the computation of the transformation matrix between them. After subdivision, node $7$ is a T-junction in the full system (see~\cref{fig:t_junction_schema}a). Within this heuristic illustration, considering that each node has only one DOF and the static condition is applied, the displacement of the node $7$ is expressed from those of node $2$ and $3$ as
\begin{equation}
\label{eq:displacement_relation}
 u_7 = 0.5 u_2 + 0.5 u_3.
\end{equation}
\cref{fig:t_junction_schema}b shows the reduced system where the T-junction node $7$ is not considered. The displacement fields between these full and reduced systems are expressed (by taking into account~\cref{eq:displacement_relation}) through the transformation matrix $\mathbf{T}$ as
\begin{equation}
\resizebox{0.77\columnwidth}{!}{$
\underbrace{
 \begin{Bmatrix}
     u_0 \\[0.3em] u_1 \\[0.3em]  \vdots \\[0.3em]  u_7
  \end{Bmatrix} }_{\mathbf{u}_f}
  = \underbrace{ \begin{pmatrix}
  1 & 0 & 0 & 0 & 0 & 0 & 0 \\[0.3em]
  0 & 1 & 0 & 0 & 0 & 0 & 0 \\[0.3em]
  0 & 0 & 1 & 0 & 0 & 0 & 0 \\[0.3em]
  0 & 0 & 0 & 1 & 0 & 0 & 0 \\[0.3em]
  0 & 0 & 0 & 0 & 1 & 0 & 0 \\[0.3em]
  0 & 0 & 0 & 0 & 0 & 1 & 0 \\[0.3em]
  0 & 0 & 0 & 0 & 0 & 0 & 1 \\[0.3em]
  0 & 0 & 0.5 & 0.5 & 0 & 0 & 0 
 \end{pmatrix} }_{\mathbf{T}}
 \underbrace{
 \begin{Bmatrix}
     u_0 \\[0.3em] u_1 \\[0.3em]  \vdots \\[0.3em]  u_6
  \end{Bmatrix}.
  }_{\mathbf{u}_r}$
  }
\end{equation}
The transformation matrix $\mathbf{T}$ for the general 3D case where each node has three DOFs is built straightforwardly from this example.

Assuming we are using a non-linear constitutive law (\eg{} hyperelastic) or a co-rotational FE formulation, the system matrix needs to be updated at each time step. Consequently, the local updating of the topology has a very limited impact on the computation. The main overhead comes from the T-junction handling, but it is somewhat compensated by the reduced dimensions of the linear system to be solved (although the reduced matrix is denser than the initial one). Our experience has shown that if we consider that about $10\%$ of the mesh elements are subdivided, approximately $20\%$ of the nodes in the resulting mesh are T-junctions. Obviously, the number of T-junctions depends strongly on the template mesh used by the refinement, and also on the fact that the elements are subdivided locally within one or several regions. 







\section{Needle-tissue interaction algorithm}
\label{sec:needle_tissue_interaction_algo}

To model the interaction of needle and tissue, we consider two different constraints: penetration (puncture) and sliding~\cite{Duriez2009}. To avoid remeshing when modeling needle-tissue interaction, we use the same constraints based approached as described in \cite{Duriez2009}. However, unlike \cite{Duriez2009}, we solve the constrained system differently. Before entering the tissue, the needle-tissue constraint is only created when the needle tip is in contact with the tissue surface. This penetration constraint is represented mathematically as $\mathscr{P}(\mathbf{x}_n,\mathbf{x}_t) \geq 0$, where $\mathbf{x}_n$, and $\mathbf{x}_t$ stand for the position of the needle, and tissue, respectively. Immediately after entering the tissue, sliding constraint between the needle and tissue is created along the needle trajectory as $\mathscr{S}(\mathbf{x}_n,\mathbf{x}_t) = 0$. When friction is considered, both $\mathscr{P}$ and $\mathscr{S}$ are nonlinear. The constraints $\mathscr{P}$ and $\mathscr{S}$ between the tissue (denoted by subscript $1$) and needle (denoted by subscript $2$) are expressed through the global coordinate system using Lagrange multipliers $\bm{\lambda}$ as follows
\begin{equation}
\label{eq:needle_tissue_interaction}
\resizebox{0.62\columnwidth}{!}{
 $\begin{pmatrix}
  \mathbf{A}_1 & \mathbf{0} & \mathbf{J}_1^T \\[0.3em]
  \mathbf{0} & \mathbf{A}_2 & \mathbf{J}_2^T \\[0.3em]
  \mathbf{J}_1 & \mathbf{J}_2 & \mathbf{0}
 \end{pmatrix} \begin{Bmatrix}
                d\mathbf{v}_1 \\[0.3em] d\mathbf{v}_2 \\[0.3em] \bm{\lambda}
               \end{Bmatrix}
               = \begin{Bmatrix}
                  \mathbf{b_1} \\[0.3em] \mathbf{b_2} \\[0.3em] \mathbf{0}
                 \end{Bmatrix}$,}
\end{equation}


\noindent where $\mathbf{A}_1$ and $\mathbf{A}_2$ are the system matrices for the soft tissue and needle, respectively;  $\mathbf{J}_1$ and $\mathbf{J}_2$ account for contraint directions between the needle and tissue. In the local coordinate system attached to the needle, constraints between needle and tissue are only expressed in two directions orthogonal to the needle shaft, resulting in a needle-tissue sliding constraint. The expressions of these constraints in the global coordinate system $\mathbf{J}_1$ and $\mathbf{J}_2$ are then built by transforming the local constraint expressions from the local coordinate system to the global one.
However, formulating the problem as~\eqref{eq:needle_tissue_interaction} leads to a non positive definite global matrix, which makes the system challenging to solve. An alternative approach, as proposed in~\cite{Hamze2016,Duriez2009}, is to solve the interaction problem in three steps: predictive motion (no interaction constraints), constraint solving, and corrective motion. However, this alternative requires the computation of matrix inverse $\mathbf{A}_1^{-1}$ and $\mathbf{A}_2^{-1}$. This approach is time consuming, especially for large systems. Unlike this method, we solve the constrained problem iteratively by using the augmented Lagrangian method~\cite{Uzawa1989}
\begin{subequations}
\label{eq:uzawa_algo}
\begin{align}
\label{eq:uzawa_algo_a}
(\mathbf{A} + \mathbf{J}^T\mathbf{W}\mathbf{J}) d\mathbf{v}^{k+1} = \mathbf{b}-\mathbf{J}^T\bm{\lambda}^k \\
\bm{\lambda}^{k+1} = \bm{\lambda}^k - \mathbf{W}\mathbf{J}d\mathbf{v}^{k+1},
\end{align}
\end{subequations}
where $\mathbf{W}$ is the penalty weight matrix with finite values. The advantage of this method is that the exact solution of the needle-tissue interaction can be obtained, as compared to the penalty method (see e.g. \cite{IvoBabuska73}), and no additional DOFs are needed, as compared to the classical Lagrange multiplier method \cite{Papadopoulos98}. A critical feature of this approach is that the system matrix in~\eqref{eq:uzawa_algo_a} is positive definite, therefore iterative solvers, such as the conjugate gradient, can be used efficiently.
%

It is worth stressing that using the augmented Lagrangian method for solving needle-tissue interaction, combined with the T-junction handling in the tissue, is straightforward. Indeed, as mentioned above, it is sufficient to solve~\eqref{eq:uzawa_algo_a} in the reduced space, and $d\mathbf{v}^{k+1}$ is easily computed from the reduced solution.

\section{Results and discussions}
\label{sec:results}


To demonstrate the efficiency of our method, we present several numerical studies. We first present the convergence of the stress error on a typical L-shaped domain.
The motivation for this test is because of its localized nature, \ie{} stress is concentrated at the corner of the domain, which mimics the localized scenario of needle insertion. To demonstrate the computational advantage of adaptive refinement over uniform refinement, we also present the computational time for this problem.
%
Then, to point out the benefits of a local mesh refinement in needle insertion simulation, we study a needle insertion scenario with friction to show the impact of local refinement on the displacement field around the needle shaft and also on the needle-tissue interaction force profile. Finally, a more complicated scenario is simulated \ie{} insertion of a needle into a liver which is undergoing breathing motion. In our simulations, the needle and soft tissue follow a linear elastic constitutive law, associated to a co-rotational FE formulation. 

\subsection{Convergence study}

\begin{figure}[!htbp]
\centering
\includegraphics[width=0.5\columnwidth]{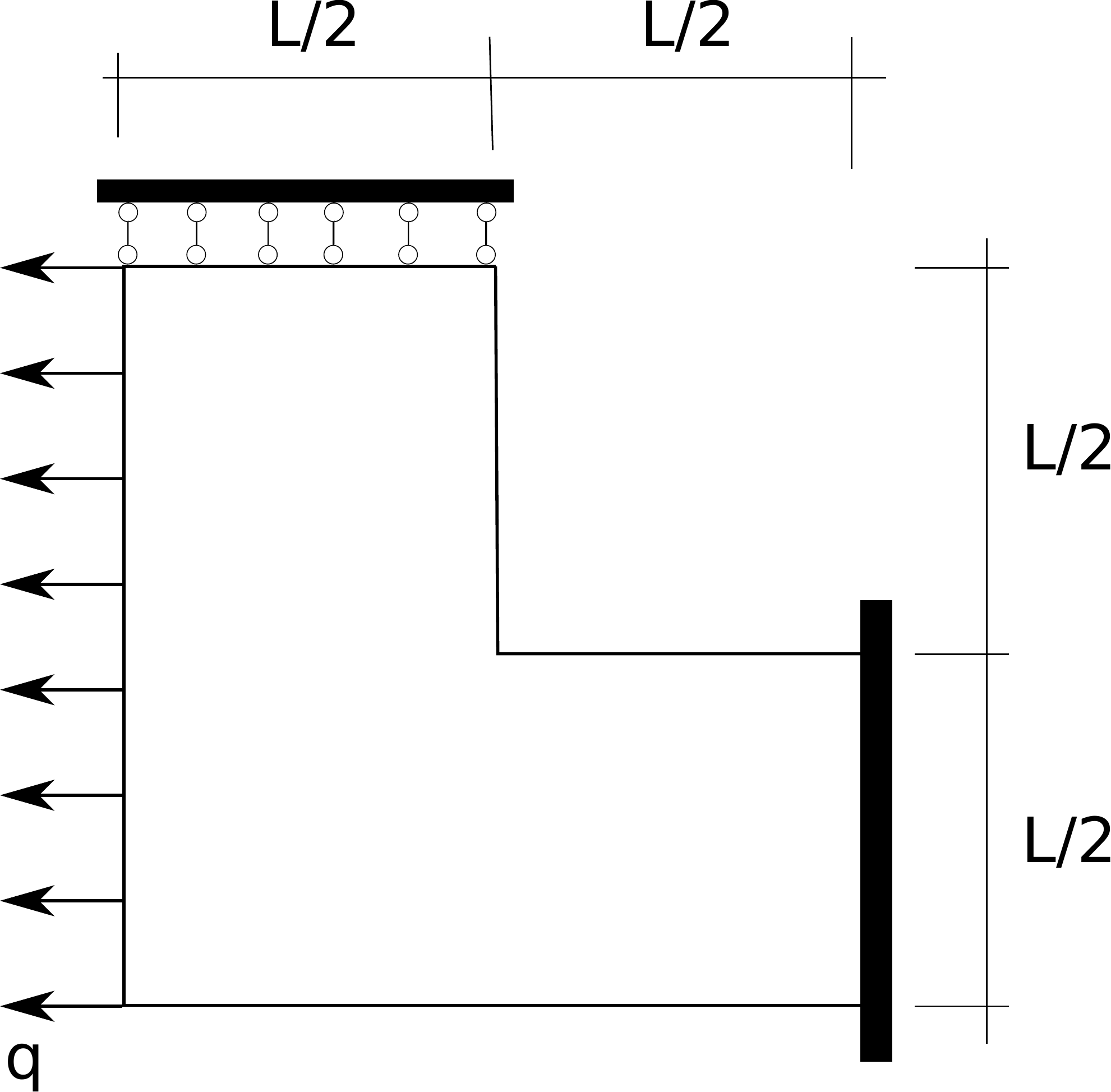}
\caption{Boundary conditions of the L-shaped domain test.}
\label{fig:L_shape}
\end{figure}

To show the advantage of error-controlled adaptive refinement scheme, as compared to the uniform mesh refinement, a convergence study is performed on a 3D L-shaped domain.
As shown in~\cref{fig:L_shape}, the L-shaped domain is clamped at the right boundary, and simply supported in the vertical direction at the top boundary. The dimension is set to $L=4$ and the thickness of the domain is $L/2 = 2$. Young's modulus, and Poisson's ratio of the tested material are $1\times 10^3$, and $0.3$, respectively. The domain is subjected to a uniformly distributed traction force on the left surface boundary.

\begin{figure}[!htbp]
\centering
\includegraphics[width=0.86\columnwidth]{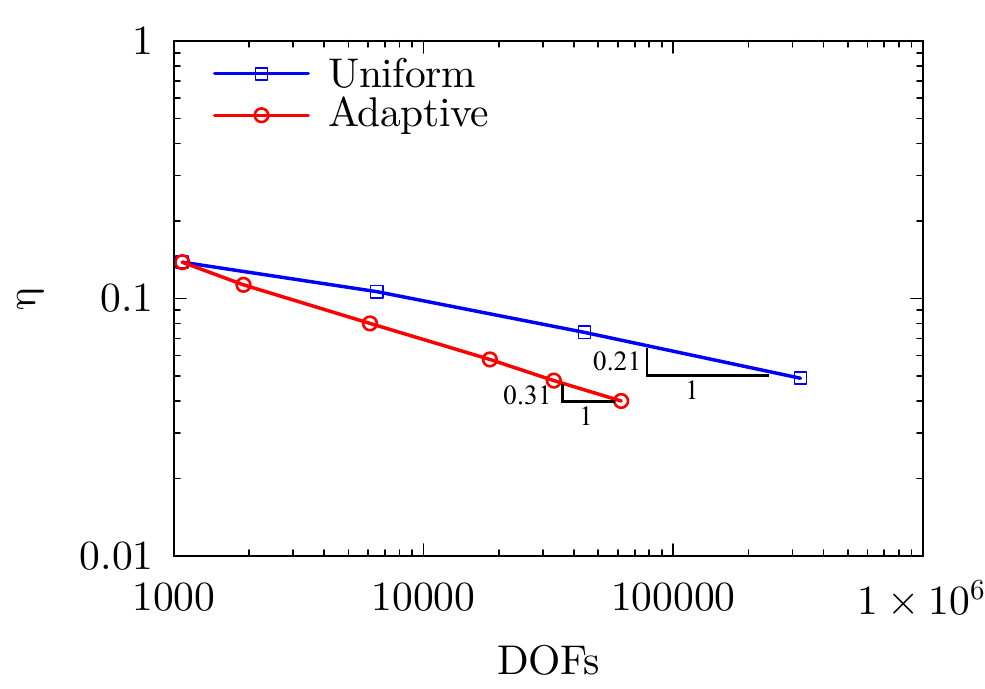}
\caption{Convergence of the relative error, comparison between uniform and adaptive refinements.}
\label{fig:convergence}
\end{figure}

Starting with the mesh having $8 \times 8 \times 4$ hexahedral elements (excluding the $4 \times 4 \times 4$ corner elements),
two types of refinements are performed. The first one, called uniform refinement, consists of subsequently subdividing every hexahedral element in to $8$ smaller elements. In the second approach, called adaptive refinement, only those elements which satisfy the marking condition \eqref{eq:max_rule} are refined (subdivided in to $8$ smaller elements).
%
%
We define the relative error $\eta$ as
\begin{equation}
\eta =  \frac{  \sqrt{ \sum_{e=1}^{N_{e}} \int_{\Omega^e} (\bm{\epsilon}^h - \bm{\epsilon}^s)^T (\bm{\sigma}^h - \bm{\sigma}^s)  \mathrm d \Omega } }
{ \sqrt{ \sum_{e=1}^{N_{e}} \int_{\Omega^e} (\bm{\epsilon}^h)^T \bm{\sigma}^h \mathrm d \Omega } }.
\end{equation}
%
%
In \cref{fig:convergence}, we show the plots of the relative error versus the number of DOFs for uniform and adaptive refinement.
We see that for the uniform refinement, the relative error $\eta$ converges with a slope of $0.21$, which corresponds to the theoretical slope of $2/9$ for singular problems in 3D. By comparison with the uniform refinement, the adaptive refinement converges with a higher slope ($0.31$). Clearly, to achieve certain expected error of the simulation, the adaptive refinement needs fewer DOFs than the uniform refinement.
%
%

%

%
To demonstrate the performance of adaptive (local) refinement in terms of computational time, as compared to uniform (full) refinement of the mesh, with the same expected relative error $\eta=8\%$, we again studied the L-shaped domain problem.
The result is reported in~\cref{tab:computation_time}. The local refinement decreases the number of DOFs by a factor $5.9$ (7473 vs. 44064) associated with a computational speed-up of $16\times$ (736.41 ms vs. 12140.8 ms).
%
%
\begin{table}[!htbp]
 \centering
 \begin{tabular}[c]{ l | c | c | c | c }
 {} &  \multirow{2}{*}{DOFs}   &    \multicolumn{2}{ c|}{Time}   &   Total \\  \cline{3-4}
 {} &  & TC & MS & time \\ \hline
 Full refinement & 44064 & x & 12140.8 & 12140.8 \\  \hline
 Local refinement & 7473 & 319.24 & 417.17 & 736.41 \\  \hline  
 \end{tabular}
 \caption{Computational time (in $10^{-3}s$). TC: Topological Changes, MS: System Matrix Solve.}
\label{tab:computation_time}
\end{table}

%

%
In view of above observations, it is a strong argument to support the employment of adaptive refinement scheme while limiting the discretization error in real-time simulations.
%


\subsection{Impact of local mesh refinement on displacement field}
\label{sec:impact_on_displacement}
\begin{figure*}[!htbp]
\centering
 \includegraphics[width=1\textwidth]{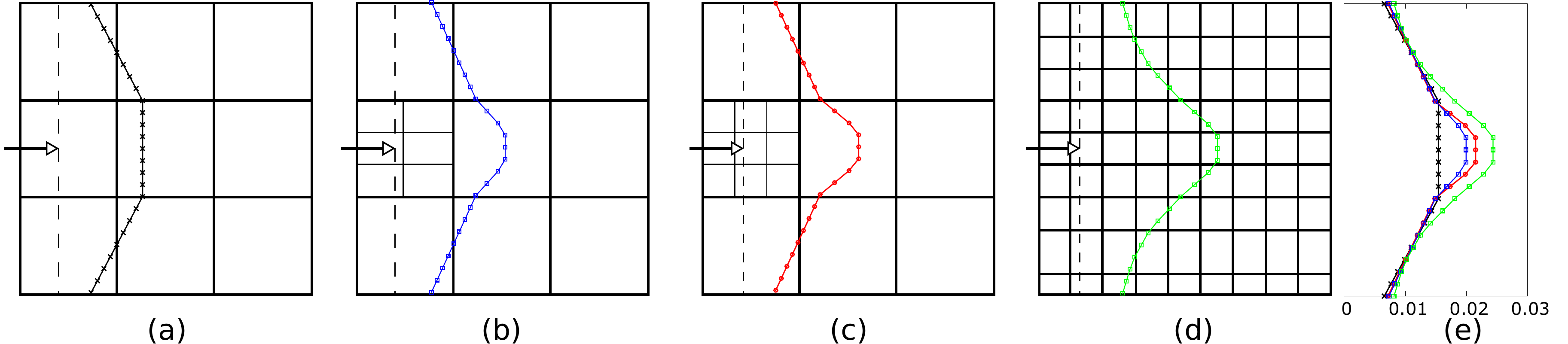}
\caption{Variation of tissue displacement resulting from friction during needle insertion, measured along a vertical line located at the needle tip. (a) The un-refined mesh; (b) adaptive refinement using an anisotropic template $2 \times 3 \times 3$; (c) adaptive refinement using an isotropic template $3 \times 3 \times 3$; and (d) full refinement.
The graph (e) shows the benefits of the anisotropic refinement.}
\label{fig:displacement_field_study}
\end{figure*}
We now present the results of a simulation of needle insertion into a homogeneous 3D tissue model. For this study, we consider the Young's modulus of  $10^8$ for the needle, and $10^3$ for the tissue, whereas the Poisson's ratio is taken as $0.4$ for both. The friction coefficient between the needle and the tissue is set to $0.9$. The displacement field due to frictional interactions with the needle, viewed from the $xy$ plane of the tissue, is shown in \cref{fig:displacement_field_study}.
It is shown (by the nonlinear variation of the displacement field in the vicinity of the needle) that when the mesh is adaptively refined near the needle shaft, the needle-tissue interaction is captured as good as in the case of full refinement (see \cref{fig:displacement_field_study}d). Indeed, closer the position is to the needle shaft, higher the obtained displacement field. Conversely, when a coarse element is used, and is not refined during the simulation, the above behavior is not reproduced within the element (\cref{fig:displacement_field_study}a). It is important to point out that the refinement using anisotropic template, as in~\cref{fig:displacement_field_study}b, is very relevant since it generates fewer DOFs than using the isotropic template (\cref{fig:displacement_field_study}c), while still catching the nonlinear displacement field.

\subsection{Impact of local mesh refinement on needle-tissue interaction}
\label{sec:needle_beam_insertion}

\begin{figure}[!htbp]
\centering
 \includegraphics[width=0.8\columnwidth]{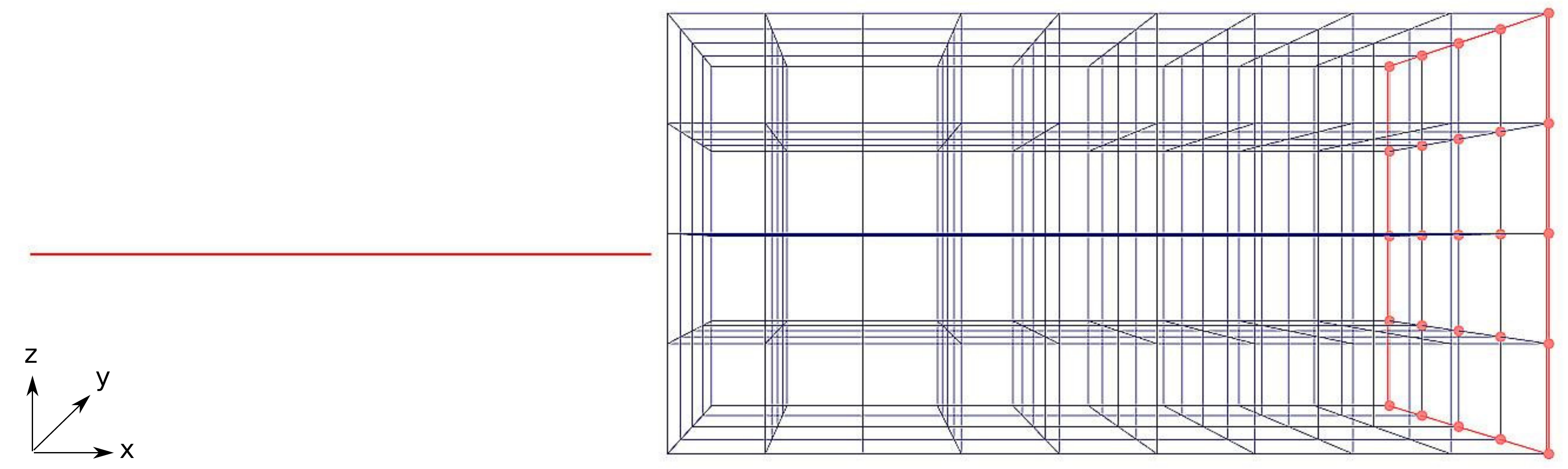}
\caption{Schematic representation of needle insertion simulation into a phantom tissue. The phantom tissue is clamped at the right surface.}
\label{fig:needle_insertion_beam_sharpen}
\end{figure}

In order to gain insight into the nonlinear behavior of the needle-tissue interaction around the needle shaft, and to exhibit the effect of the adaptive refinement on a trade-off between
computational time and precision, a needle insertion simulation into a phantom tissue test is carried out, see~\cref{fig:needle_insertion_beam_sharpen}. For this study, we consider the Young's modulus of  $50$~MPa for the needle, and $10$~MPa for the tissue, whereas the Poisson's ratio is taken as $0.4$ for the tissue, and $0.3$ for the needle, respectively. Again, a linear elastic model based on corotational formulation is employed for the needle as well as the tissue. The dimension of the tissue is $4 \times 2 \times 2$~cm. The needle length and radius are of $3.2$~cm and $0.1$~cm, respectively.

Three meshing schemes are employed: a coarse mesh with resolution $10 \times 5 \times 5$ nodes, a fine mesh with resolution $20 \times 10 \times 10$ nodes and an adaptive mesh (starting with the coarse mesh $10 \times 5 \times 5$ nodes and adaptively refining the mesh during the simulation). Within the adaptive meshing scheme, the mesh refinement is again piloted by the error estimate described in~\cref{sec:Error_estimate_and_adaptive_refinement}.

To investigate the sensitivity of the needle-tissue interaction parameters (frictional coefficient and puncture strength) on the resulting mesh adaptation, and thus on computational output, two scenarios are studied. The first concerns varying the puncture strength $\sigma_{n0}$ parameter ($0$~N, $10$~N, and $20$~N), while keeping the same frictional coefficient $\mu = 0.5$ between the tissue and the needle shaft. The second is dedicated to study the influence of the frictional coefficient by setting it to $0.1$, $0.3$ and $0.5$, while keeping the puncture strength unchanged ($\sigma_{n0} = 10$~N). Within these two scenarios, the frictional coefficient on the tissue surface is set to $0.8$.

\begin{figure*}[!htbp]
 \centering
      \begin{subfigure}[b]{0.32\textwidth}
              \centering
              \includegraphics[width=1\columnwidth]{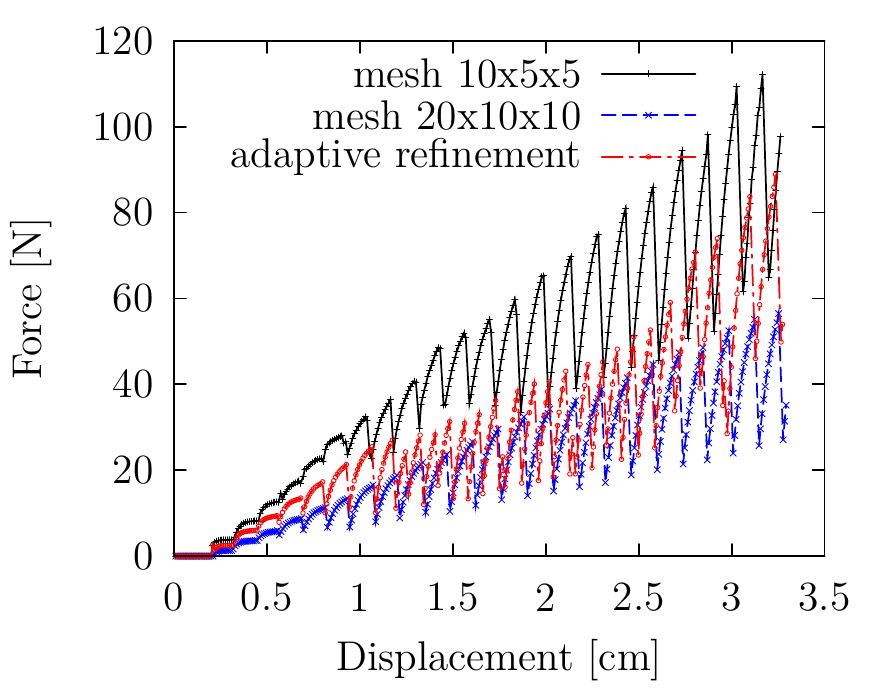}
              \caption{$\lambda_{n0}$ = 0~N}
              \label{fig:penetration_force_00}
      \end{subfigure}
        ~ 
      \begin{subfigure}[b]{0.32\textwidth}
              \centering
	      \includegraphics[width=1\columnwidth]{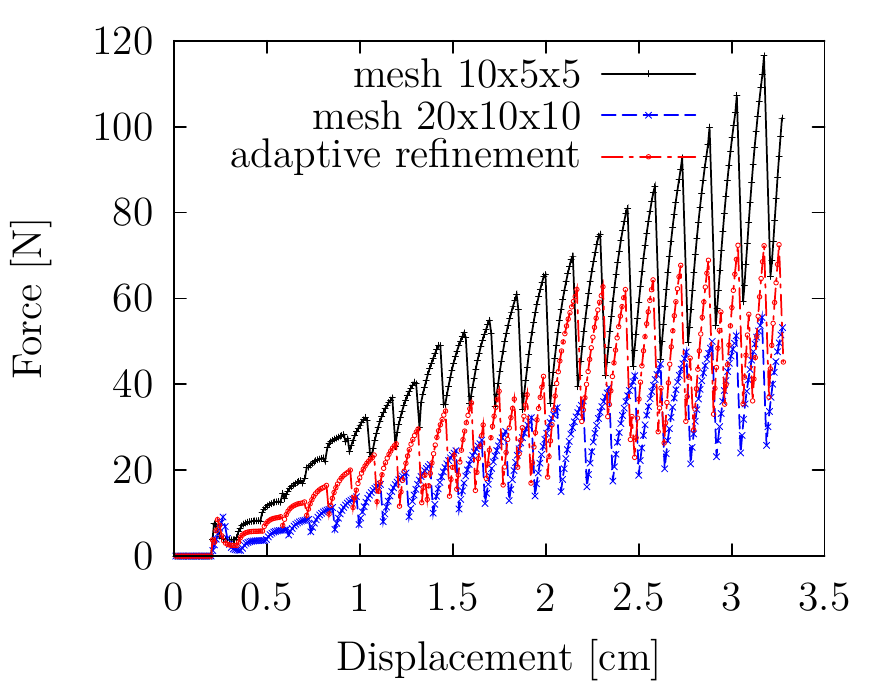}
              \caption{$\lambda_{n0}$ = 10~N}
              \label{fig:penetration_force_10}
      \end{subfigure}%
       ~ 
      \begin{subfigure}[b]{0.32\textwidth}
              \centering
	      \includegraphics[width=1\columnwidth]{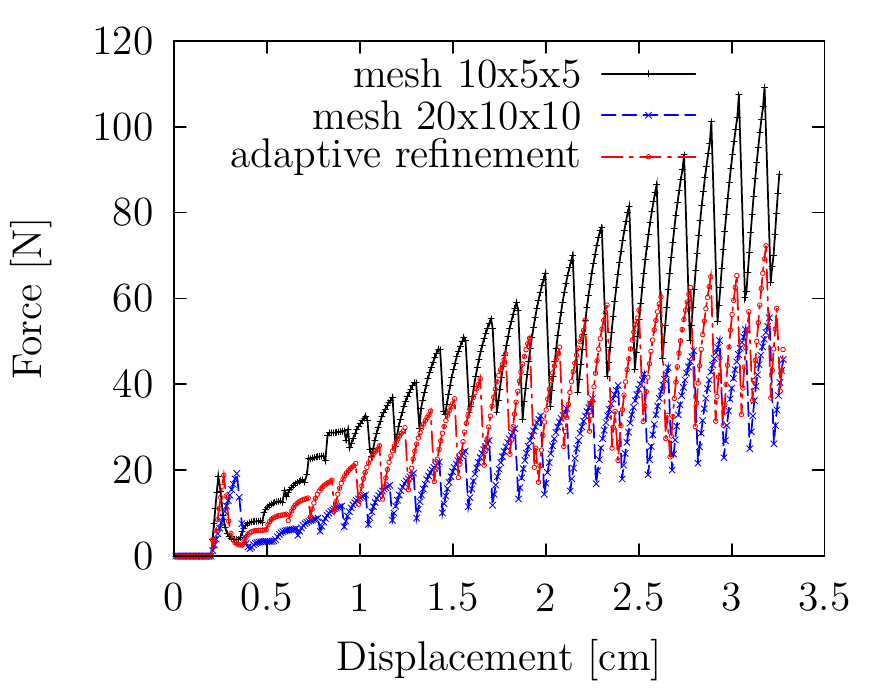}
              \caption{$\lambda_{n0}$ = 20~N}
              \label{fig:penetration_force_20}
      \end{subfigure}%
\caption{Comparison of needle-tissue interaction forces along the needle shaft within two cases: without refinement (with different mesh resolutions) and with adaptive refinement. The penetration strength is varied while keeping the same frictional coefficient $\mu=0.5$ between the needle shaft and the soft tissue.}\label{fig:force_displacement_needle}
\end{figure*}

\begin{figure*}[!htbp]
 \centering
       \begin{subfigure}[b]{0.32\textwidth}
              \centering
	      \includegraphics[width=1\columnwidth]{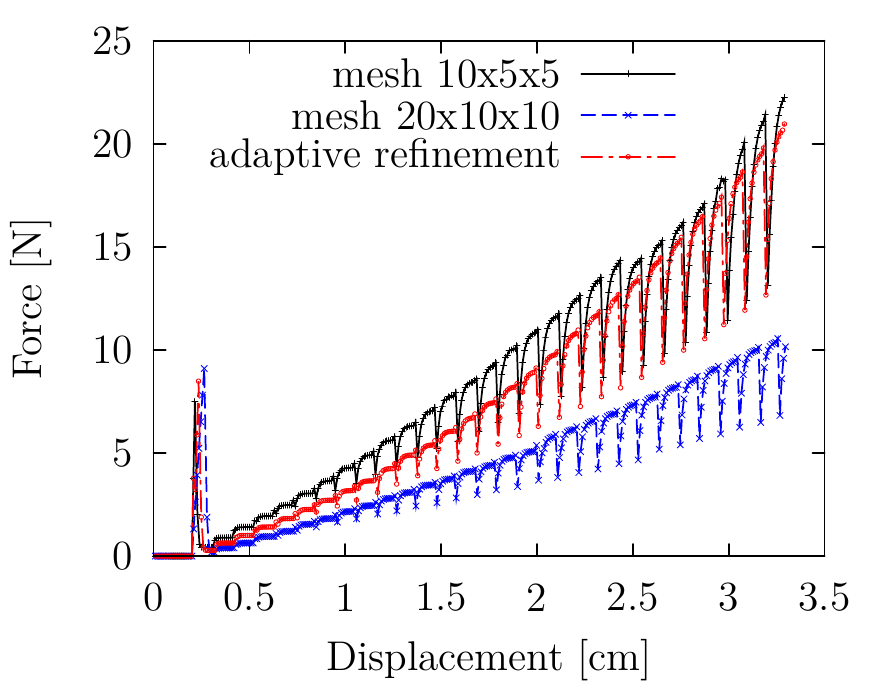}
              \caption{$\mu = 0.1$}
              \label{fig:penetration_force_10_friction_0.1}
      \end{subfigure}%
             ~ 
      \begin{subfigure}[b]{0.32\textwidth}
              \centering
	      \includegraphics[width=1\columnwidth]{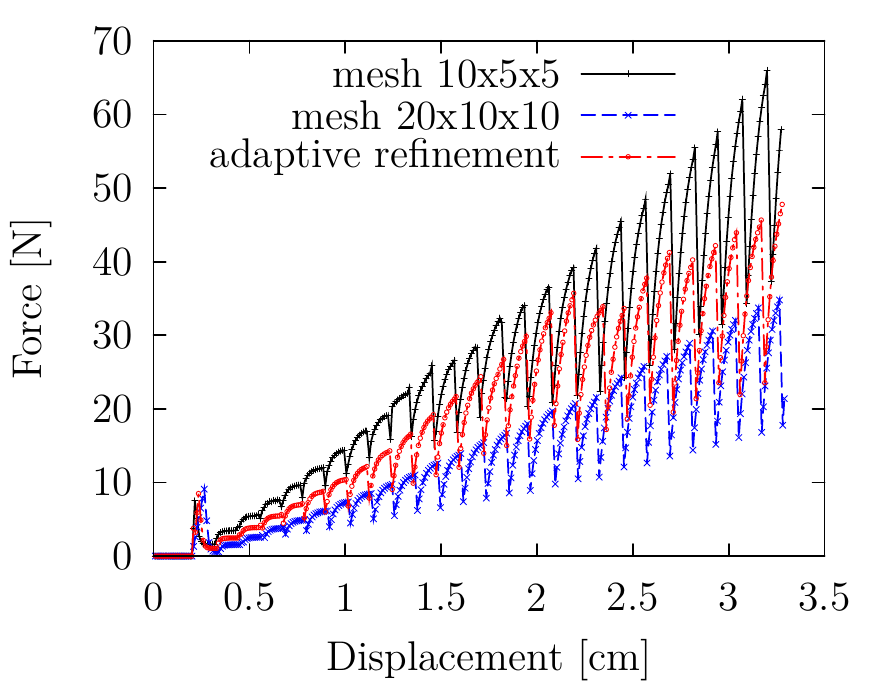}
              \caption{$\mu = 0.3$}
              \label{fig:penetration_force_10_friction_0.3}
      \end{subfigure}%
       ~ 
      \begin{subfigure}[b]{0.32\textwidth}
              \centering
	      \includegraphics[width=1\columnwidth]{friction_05_penetration_force_10_force_displacement.pdf}
              \caption{$\mu = 0.5$}
              \label{fig:penetration_force_10_friction_0.5}
      \end{subfigure}%
\caption{Comparison of needle-tissue interaction forces along the needle shaft within two cases: without refinement (with different mesh resolutions) and with adaptive refinement. The frictional coefficient between the needle shaft is varied while keeping the same penetration strength (10~N) at the tissue surface.}\label{fig:force_displacement_needle_friction}
\end{figure*}

\cref{fig:force_displacement_needle} shows the plots of the integrated interaction force along the needle shaft versus the displacement of the needle tip for the first scenario. The second scenario is depicted in~\cref{fig:force_displacement_needle_friction}.

\begin{figure*}[!htbp]
 \centering
      \begin{subfigure}[b]{0.4\textwidth}
              \centering
	      \includegraphics[width=1\columnwidth]{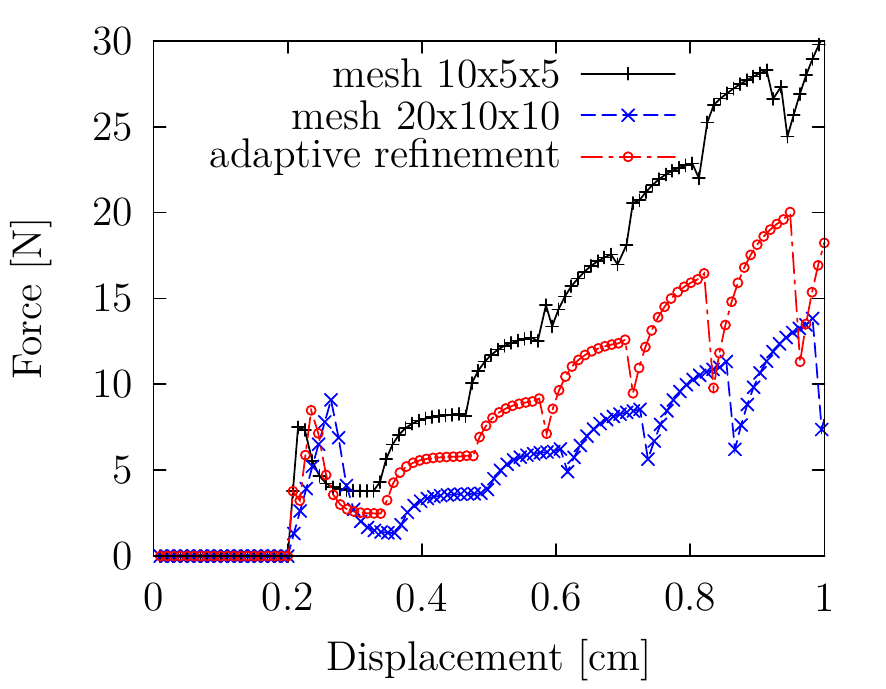}
	      \caption{}
	      \label{fig:force_displacement_zoom_a}
      \end{subfigure}
        ~ 
      \begin{subfigure}[b]{0.4\textwidth}
              \centering
	      \includegraphics[width=1\columnwidth]{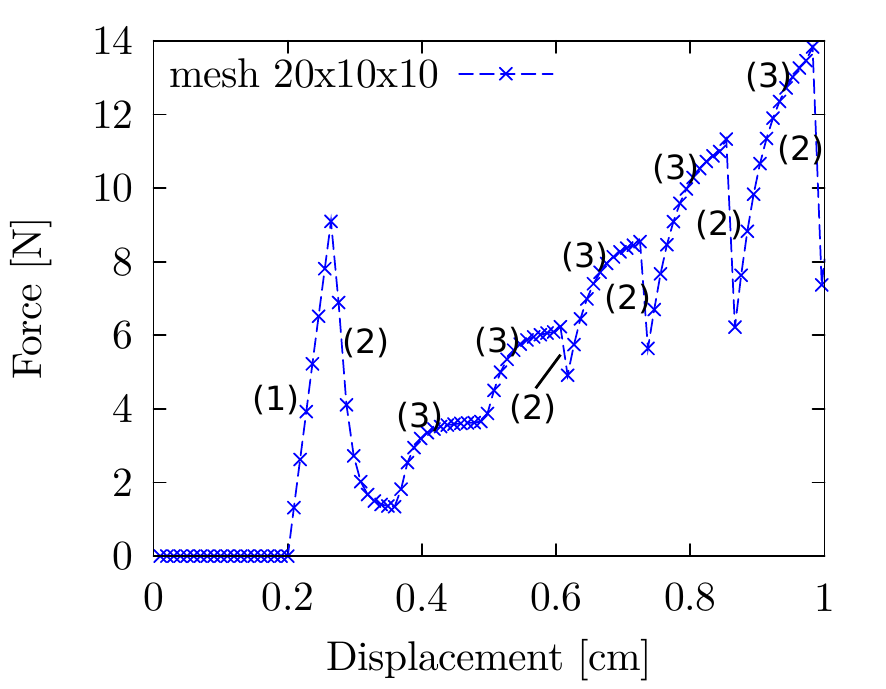}
	      \caption{}
	      \label{fig:force_displacement_zoom_mesh20x10x10_behavior}
      \end{subfigure}%
\caption{(\subref{fig:force_displacement_zoom_a}) Puncture, cutting and relaxation behaviors are shown by a zoom-in plot from \cref{fig:penetration_force_10}. (\subref{fig:force_displacement_zoom_mesh20x10x10_behavior}) The typical behavior is shown in phases. Phase (1): The needle is puncturing the tissue surface. Phase (2): Just after the penetration event, the relaxation occurs. Phase (3): The interaction force increases due to the fact that frictional force increases with insertion distance. When the needle tip has cut the tissue to advance forward, the relaxation occurs again.}
\label{fig:force_displacement_zoom}
\end{figure*}

It shows that when the contact force between the needle tip and the tissue surface is higher than the tissue puncture strength, the needle penetrates into the tissue. Right after this penetration event, a relaxation phase can be observed that induces a decreasing force at the needle bases. Thereafter, it is observed that as the needle moves forward, the interaction force increases due to the increasing frictional force along the needle shaft (which is directly proportional to the insertion distance). Only when the contact force at the needle tip is greater than the cutting strength of the soft tissue, the needle cuts the tissue and continues going ahead. Immediately after this cutting action, the relaxation phase is observed anew. This behavior is periodically observed during the needle insertion. These observations are clearly shown in~\cref{fig:force_displacement_zoom_a}, which is obtained by zooming in~\cref{fig:penetration_force_10}. A typical behavior with distinguished phases is presented in~\cref{fig:force_displacement_zoom_mesh20x10x10_behavior}.

\begin{figure}[!htbp]
\centering
 \includegraphics[width=0.75\columnwidth]{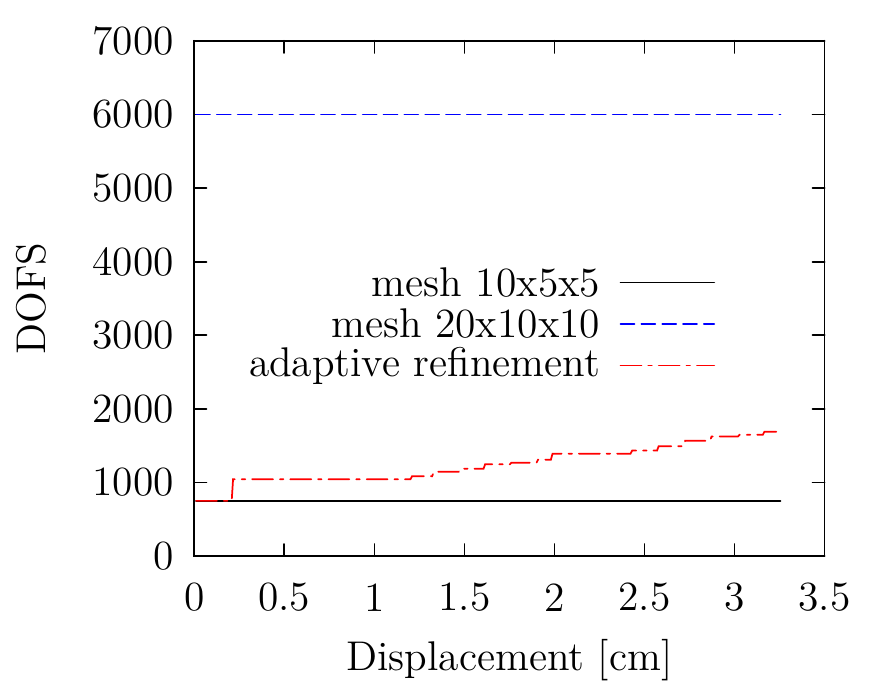}
\caption{Number of DOFs during the needle insertion of the simulation in~\cref{fig:penetration_force_10}.}
\label{fig:number_dofs}
\end{figure}

It is observed that under mesh refinement the resulting global behavior of needle-tissue system is less stiff. This is explained by the fact that beneath mesh refinement, a greater displacement field is obtained, which results from the needle-tissue interaction (as also observed in~\cref{sec:impact_on_displacement}). It is also shown that using the adaptive refinement scheme, the interaction needle-tissue behavior is close to those when using a fine mesh (see~\cref{fig:force_displacement_needle}). As can be seen in~\cref{fig:number_dofs}, an interesting observation is that the number of DOFs in the adaptive refinement simulation is significantly fewer than that of using the fine mesh. This obviously results in an important gain in terms of computational time. Indeed, the simulation using the adaptive refinement mesh runs at nearly $45$~FPS compared to $4$~FPS of that using the uniform fine mesh.

\begin{figure*}[!htbp]
 \centering
       \begin{subfigure}[b]{0.3\textwidth}
              \centering
	      \includegraphics[width=1\columnwidth]{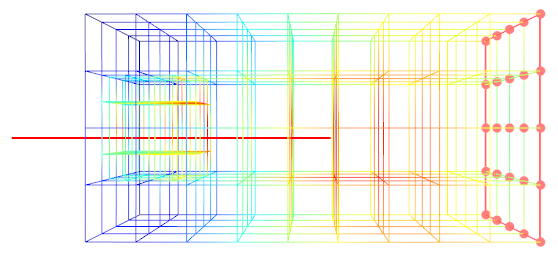}
              \caption{$\mu = 0.1$}
              \label{fig:friction_01}
      \end{subfigure}%
             ~ 
      \begin{subfigure}[b]{0.3\textwidth}
              \centering
	      \includegraphics[width=1\columnwidth]{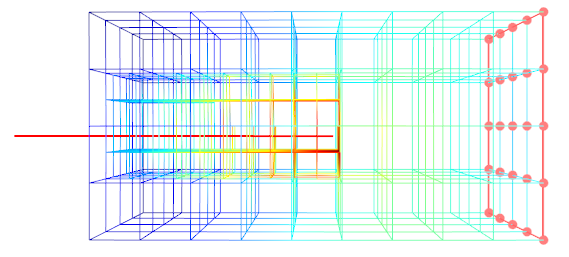}
              \caption{$\mu = 0.3$}
              \label{fig:friction_03}
      \end{subfigure}%
       ~ 
      \begin{subfigure}[b]{0.3\textwidth}
              \centering
	      \includegraphics[width=1\columnwidth]{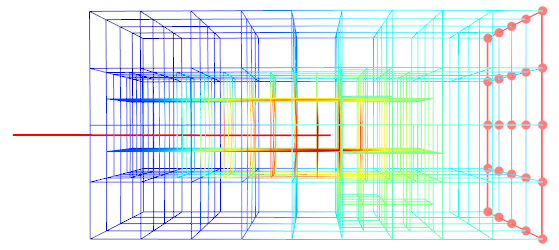}
              \caption{$\mu = 0.5$}
              \label{fig:friction_05}
      \end{subfigure}%
\caption{Refinement patterns (colored by stress level) of the adaptive scenario while using different frictional coefficients.}\label{fig:force_displacement_needle_friction_patterns}
\end{figure*}

It is seen from~\cref{fig:force_displacement_needle_friction} that the smaller the frictional coefficient is, the more the behavior of the adaptive refinement scheme differs from that of the simulation using the fine mesh. It aligns nicely with the fact that smaller friction force does not lead to mesh refinement around the needle shaft. Indeed, as seen in~\cref{fig:force_displacement_needle_friction_patterns}, the refinement in the case of frictional coefficient $\mu = 0.1$ is mostly due to the penetration force at the tissue surface.

\subsection{Application to liver}

The method proposed in this paper is now applied to a liver model undergoing breathing motion, to mimic a typical case of radio-frequency ablation of a tumor. The same Young's modulus and Poisson's ratio for the needle and the tissue as in~\cref{sec:needle_beam_insertion} are employed. The frictional coefficient is set to $0.5$ when the needle is inserted and to $0.1$ when it is pulled back. The puncture force at the tissue surface is set to $10$~N. Induced by error estimate (\cref{eq:max_rule}), the needle insertion and constraints applied to the liver lead to refinements in different regions. The initial mesh has 1179 DOFs. When the needle advances into the liver, combining with the motion of the liver due to breathing effect, the mesh is progressively refined to accurately take in to account the interaction of needle and liver. The maximum number of DOFs when the needle is completely inserted into the liver is 2961. And, when the needle is steadily pulled back, the mesh is then progressively coarsened until the needle is completely outside of the liver. Thereafter, the refinement process is now only due to the movement of the liver by breathing effect and imposed boundary conditions. The number of DOFs at this stage is 1509. By applying this adaptive refinement\slash coarsening procedure, it is not only guaranteed that the discretization error is fully controlled, but the computational cost is also kept as small as possible. Indeed, without adaptive remeshing procedure applied on the initial mesh, the simulation runs at 35 FPS while the discretization error is $12\%$, whereas when the adaptive refinement is performed, it runs at 22 FPS while decreasing the discretization error to $6\%$. Note that these frame rates result not only from computational resolutions of needle, tissue and their interactions but also from their visualization cost.
%
\begin{figure*}[!htbp]
 \centering
 \def\svgwidth{0.78\textwidth}
 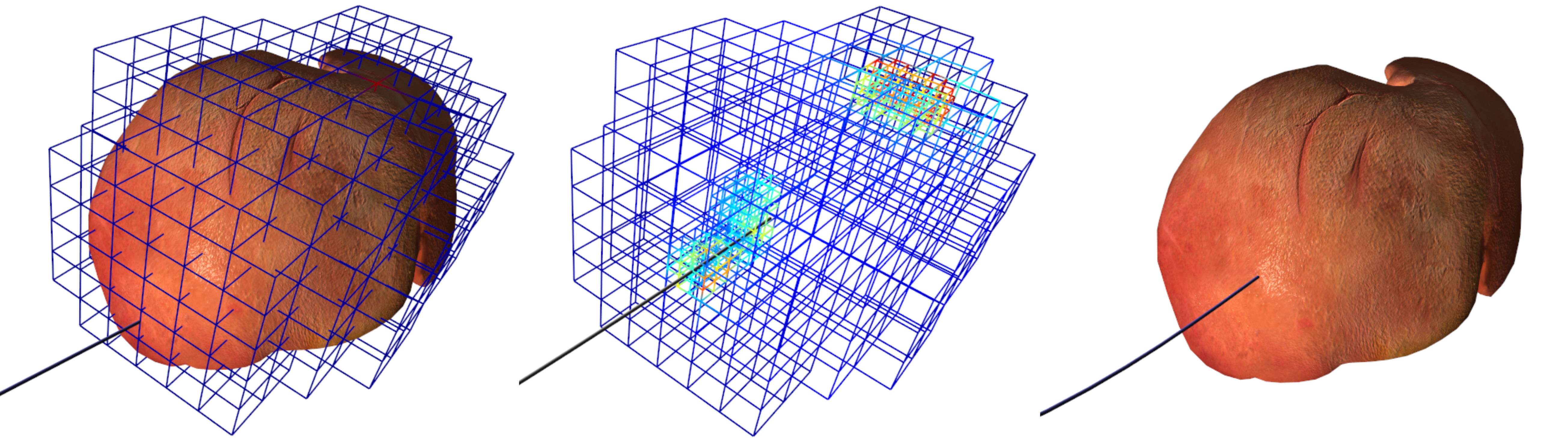
 \caption{(a) Simulation of needle insertion in a liver; (b) using dynamic mesh refinement scheme driven by error estimate; (c) visual depiction. The simulation runs at 22 Hz (on a 4GHz processor).}
\label{fig:liver_application}
\end{figure*}
%

\begin{figure}[!htbp]
\centering
 \includegraphics[width=0.91\columnwidth]{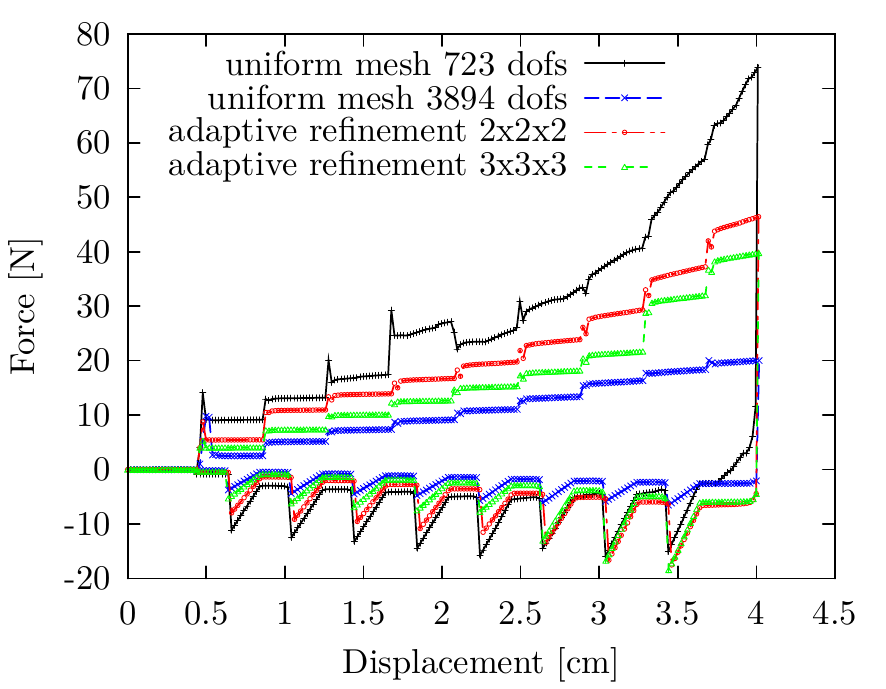}
\caption{Needle-liver phantom interaction force during needle insertion and pullback. The interaction force varies due to advancing friction and tissue cutting strength, globally increases (with positive values) during the insertion stage. At $4$~cm of the needle-tip displacement, the needle is retracted and the interaction force changes the direction, varies due to retrograding friction, globally decreases and gets zero when the needle is completely pulled out. These form a loop profile for force-displacement curve resulting from needle insertion-retraction.}
\label{fig:liver_force_displacement}
\end{figure}

In order to investigate the benefits of the adaptive refinement scheme when the needle is inserted and retracted into the liver phantom, tests with uniform and adaptive refinement schemes are carried out. Within the uniform refinement case, a coarse mesh with $723$ DOFs and a fine mesh with $3894$ DOFs are used for the liver discretization. Whereas upon the adaptive refinement scenario, the simulations start with the coarse mesh $723$ DOFs and is adaptively refined by two schemata: each marked element is refined into \begin{inparaenum}[(i)] \item  $2 \times 2 \times 2$ elements, and \item into  $3 \times 3 \times 3$ elements\end{inparaenum}. The integrated needle-tissue interaction force along the needle shaft is plotted versus the needle tip displacement when the needle is inserted and pulled back, see~\cref{fig:liver_force_displacement}.

It is observed that when the needle is outside the tissue, there is no interaction force between them. This is also detected when the needle is completely retracted from the tissue. It is clearly shown that the needle-tissue interaction depends strongly on the mesh used (especially if the mesh is coarse). The mesh influence reveals a stronger effect on the insertion stage than the pullback one. This is fully understood by the fact that the frictional coefficient between the needle shaft and the tissue is more important during the insertion steps than that during the pullback ones ($0.5$ versus $0.1$, respectively). Using the coarse mesh, the puncture force at the tissue surface is not well captured compared to the case where the fine mesh or adaptive refined mesh is employed, see~\cref{fig:liver_force_displacement}. Under mesh refinement around the needle shaft guided by error estimate, the needle-tissue interaction converges to the solution of the fine mesh. However, the maximum number of DOFs when using the adaptive refinement schemes $2 \times 2 \times 2$ and $3 \times 3 \times 3$ is $1071$ and $2193$ respectively. Therefore, as observed in~\cref{sec:needle_beam_insertion}, using the adaptive refinement scheme results in significantly fewer DOFs compared to the employment of uniform fine mesh. All of these lead to the conclusion that even starting with the coarse mesh but employing the adaptive refinement scheme, the needle-tissue interaction can be simulated more precisely compared to the coarse mesh, and with significantly lower computational cost compared to the uniform fine mesh.


\section{Conclusion and perspectives}
\label{sec:conclusion}
This paper contributes a structured approach to answering the important but rarely tackled question of accuracy in surgical simulation.

The novelty of our paper is to drive local adaptive mesh refinement during needle insertion by a robust a posteriori estimate of the discretization error.

This can be seen as a first step to control the error associated with acceleration methods in needle surgical simulations, and to separate the modeling (are we solving the right problem?) and discretization error (are we solving the problem right?). 

Verification of the discrete scheme is guaranteed in our approach because the a posteriori estimate asymptotically converges to the exact error. As we use an implicit approach, we also control the error on the equilibrium equations. As such, assuming a proper material model and kinematics for the problem, we can guarantee the accuracy of the solution. This is not the case in explicit time stepping approaches.

Validation of the approach is not considered. We focus here only on one source of error (discretization). Whilst this is a limitation, we do believe that quantifying discretization errors separately to modeling errors is necessary to devise accurate and clinically-usable surgical simulators and to better understand the resulting simulation results.
 
A natural direction for research, building on recent work on data-driven simulations is to devise error-controlled approaches able to learn from data as it is acquired during the simulation \cite{seiler2012enriching,seiler2014data,seiler2014efficient}. In such a paradigm, the model would adapt to the real situation, as opposed to being driven by a continuous indirect comparison, as is the case in this work, with an unknown exact solution. In turn, such a data-driven approach would facilitate patient-specific simulations, which were not considered here. We are currently investigating such directions through Bayesian inference for parameter identification and model selection~\cite{rappel2016bayesian} and uncertainty quantification approaches~\cite{Hauseux2017}.

\section*{Acknowledgements}
The first and last author are supported by the Fellowship of the last author as part of the University of Strasbourg Institute for Advanced Study (BPC 14/Arc 10138).

St\'ephane P.A. Bordas and Satyendra Tomar also thank partial funding for their time provided by the European Research Council Starting Independent Research Grant (ERC Stg grant agreement No. 279578) RealTCut towards real time multiscale simulation of cutting in non-linear materials with applications to surgical simulation and computer guided surgery".

Inria thanks for funding of European project RASimAs (FP7 ICT-2013.5.2 No610425). 

St\'ephane P.A. Bordas is grateful for many helpful discussions with Dr. Pierre Kerfriden, Prof. Karol Miller, Prof. Christian Duriez, Prof. Michel Audette and Mr. Diyako Ghaffari.
\ifCLASSOPTIONcaptionsoff
  \newpage
\fi



\bibliographystyle{IEEEtran}
\bibliography{IEEEabrv,refs}

\begin{thebibliography}{10}
\providecommand{\url}[1]{#1}
\csname url@samestyle\endcsname
\providecommand{\newblock}{\relax}
\providecommand{\bibinfo}[2]{#2}
\providecommand{\BIBentrySTDinterwordspacing}{\spaceskip=0pt\relax}
\providecommand{\BIBentryALTinterwordstretchfactor}{4}
\providecommand{\BIBentryALTinterwordspacing}{\spaceskip=\fontdimen2\font plus
\BIBentryALTinterwordstretchfactor\fontdimen3\font minus
  \fontdimen4\font\relax}
\providecommand{\BIBforeignlanguage}[2]{{%
\expandafter\ifx\csname l@#1\endcsname\relax
\typeout{** WARNING: IEEEtran.bst: No hyphenation pattern has been}%
\typeout{** loaded for the language `#1'. Using the pattern for}%
\typeout{** the default language instead.}%
\else
\language=\csname l@#1\endcsname
\fi
#2}}
\providecommand{\BIBdecl}{\relax}
\BIBdecl

\bibitem{nealen2006physically}
A.~Nealen \emph{et~al.}, ``Physically based deformable models in computer
  graphics,'' in \emph{Computer graphics forum}, vol.~25, no.~4, 2006, pp.
  809--836.

\bibitem{wang2015linear}
Y.~Wang \emph{et~al.}, ``Linear subspace design for real-time shape
  deformation,'' \emph{ACM Transactions on Graphics (TOG)}, vol.~34, no.~4,
  p.~57, 2015.

\bibitem{Courtecuisse2014}
H.~Courtecuisse \emph{et~al.}, ``Real-time simulation of contact and cutting of
  heterogeneous soft-tissues,'' \emph{Medical Image Analysis}, vol.~18, no.~2,
  pp. 394 -- 410, 2014.

\bibitem{zienkiewicz1977finite}
O.~Zienkiewicz \emph{et~al.}, \emph{The finite element method: Its basis and
  fundamentals}.\hskip 1em plus 0.5em minus 0.4em\relax Elsevier, 2013, vol.~1.

\bibitem{nguyen2008meshless}
V.~P. Nguyen \emph{et~al.}, ``Meshless methods: a review and computer
  implementation aspects,'' \emph{Mathematics and computers in simulation},
  vol.~79, no.~3, pp. 763--813, 2008.

\bibitem{wu2001adaptive}
X.~Wu \emph{et~al.}, ``Adaptive nonlinear finite elements for deformable body
  simulation using dynamic progressive meshes,'' in \emph{Computer Graphics
  Forum}, vol.~20, no.~3, 2001, pp. 349--358.

\bibitem{muller2002stable}
M.~M{\"u}ller \emph{et~al.}, ``Stable real-time deformations,'' in
  \emph{Proceedings of the 2002 ACM SIGGRAPH/Eurographics symposium on Computer
  animation}, 2002, pp. 49--54.

\bibitem{ainsworth2011posteriori}
M.~Ainsworth and J.~T. Oden, \emph{A posteriori error estimation in finite
  element analysis}.\hskip 1em plus 0.5em minus 0.4em\relax John Wiley \& Sons,
  2011, vol.~37.

\bibitem{seiler2011robust}
M.~Seiler \emph{et~al.}, ``Robust interactive cutting based on an adaptive
  octree simulation mesh,'' \emph{The Visual Computer}, vol.~27, no. 6-8, pp.
  519--529, 2011.

\bibitem{Verfuerth-13-Apost}
R.~Verf{\"u}rth, \emph{A posteriori error estimation techniques for finite
  element methods}, ser. Numerical Mathematics and Scientific
  Computation.\hskip 1em plus 0.5em minus 0.4em\relax Oxford University Press,
  Oxford, 2013.

\bibitem{zienkiewicz1992superconvergent}
O.~Zienkiewicz and J.~Zhu, ``The superconvergent patch recovery (spr) and
  adaptive finite element refinement,'' \emph{Computer Methods in Applied
  Mechanics and Engineering}, vol. 101, no.~1, pp. 207--224, 1992.

\bibitem{CarstensenB-02-Apost}
C.~Carstensen and S.~Bartels, ``Each averaging technique yields reliable a
  posteriori error control in {FEM} on unstructured grids. {I}. {L}ow order
  conforming, nonconforming, and mixed {FEM},'' \emph{Math. Comp.}, vol.~71,
  no. 239, pp. 945--969 (electronic), 2002.

\bibitem{BartelsC-02-Apost}
S.~Bartels and C.~Carstensen, ``Each averaging technique yields reliable a
  posteriori error control in {FEM} on unstructured grids. {II}. {H}igher order
  {FEM},'' \emph{Math. Comp.}, vol.~71, no. 239, pp. 971--994 (electronic),
  2002.

\bibitem{BabuskaR-78-Apost}
I.~Babu{\v{s}}ka and W.~C. Rheinboldt, ``Error estimates for adaptive finite
  element computations,'' \emph{SIAM J. Numer. Anal.}, vol.~15, no.~4, pp.
  736--754, 1978.

\bibitem{Hamze2016}
N.~Hamz{\'{e}} \emph{et~al.}, ``{Preoperative trajectory planning for
  percutaneous procedures in deformable environments},'' \emph{Computerized
  Medical Imaging and Graphics}, vol.~47, pp. 16--28, 2016.

\bibitem{Abolhassani07}
N.~Abolhassani \emph{et~al.}, ``Needle insertion into soft tissue: A survey,''
  \emph{Medical Engineering and Physics}, vol.~29, pp. 413--431, 2007.

\bibitem{Duriez2009}
C.~Duriez \emph{et~al.}, ``{Interactive simulation of flexible needle
  insertions based on constraint models},'' in \emph{Lecture Notes in Computer
  Science}, vol. 5762, no. PART 2, 2009, pp. 291--299.

\bibitem{Misra2010}
S.~Misra \emph{et~al.}, ``Mechanics of flexible needles robotically steered
  through soft tissue,'' \emph{Int. J. Rob. Res.}, vol.~29, no.~13, pp.
  1640--1660, Nov. 2010.

\bibitem{wu2014physically}
J.~Wu \emph{et~al.}, ``Physically-based simulation of cuts in deformable
  bodies: A survey.'' in \emph{Eurographics (State of the Art Reports)}, 2014,
  pp. 1--19.

\bibitem{dick2011hexahedral}
C.~Dick \emph{et~al.}, ``A hexahedral multigrid approach for simulating cuts in
  deformable objects,'' \emph{IEEE Transactions on Visualization and Computer
  Graphics}, vol.~17, no.~11, pp. 1663--1675, 2011.

\bibitem{Liu201443}
G.~R. Liu and S.~S. Quek, ``{Chapter 3 - Fundamentals for Finite Element
  Method},'' in \emph{The Finite Element Method (Second Edition)}, second
  edition~ed., G.~R. Liu \emph{et~al.}, Eds.\hskip 1em plus 0.5em minus
  0.4em\relax Oxford: Butterworth-Heinemann, 2014, pp. 43--79.

\bibitem{zienkiewicz2000finite}
O.~Zienkiewicz and R.~Taylor, \emph{The Finite Element Method: Solid
  mechanics}, ser. Referex collection.Mec{\'a}nica y materiales.\hskip 1em plus
  0.5em minus 0.4em\relax Butterworth-Heinemann, 2000.

\bibitem{Felippa2005}
C.~Felippa and B.~Haugen, ``A unified formulation of small-strain corotational
  finite elements: I. theory,'' \emph{Computer Methods in Applied Mechanics and
  Engineering}, vol. 194, no. 21–24, pp. 2285 -- 2335, 2005.

\bibitem{Pinelli2010}
A.~Pinelli \emph{et~al.}, ``Immersed boundary method for generalised finite
  volume and finite difference navier-stokes solvers,'' \emph{Journal of
  Computational Physics}, vol. 229, no.~24, pp. 9073--9091, 2010.

\bibitem{Baraff1998}
D.~Baraff and A.~Witkin, ``Large steps in cloth simulation,'' in
  \emph{Proceedings of SIGGRAPH}, 1998, pp. 43--54.

\bibitem{Kwak2003531}
D.-Y. Kwak and Y.-T. Im, ``Hexahedral mesh generation for remeshing in
  three-dimensional metal forming analyses,'' \emph{Journal of Materials
  Processing Technology}, vol. 138, no. 1–3, pp. 531 -- 537, 2003.

\bibitem{Koschier2014}
D.~Koschier \emph{et~al.}, ``{Adaptive Tetrahedral Meshes for Brittle Fracture
  Simulation},'' in \emph{Symposium on Computer Animation}, V.~Koltun and
  E.~Sifakis, Eds.\hskip 1em plus 0.5em minus 0.4em\relax The Eurographics
  Association, 2014.

\bibitem{Burkhart2010}
D.~Burkhart \emph{et~al.}, ``{Adaptive and Feature-Preserving Subdivision for
  High-Quality Tetrahedral Meshes},'' \emph{{Computer Graphics Forum}},
  vol.~{29}, no.~{1}, pp. {117--127}, {2010}.

\bibitem{Paulus2015}
C.~Paulus \emph{et~al.}, ``\BIBforeignlanguage{English}{Virtual cutting of
  deformable objects based on efficient topological operations},''
  \emph{\BIBforeignlanguage{English}{The Visual Computer}}, vol.~31, no. 6-8,
  pp. 831--841, 2015.

\bibitem{Sifakis2007}
E.~Sifakis \emph{et~al.}, ``Hybrid simulation of deformable solids,'' in
  \emph{Proc. Symposium on Computer Animation}, 2007, pp. 81--90.

\bibitem{Uzawa1989}
H.~Uzawa and K.~J. Arrow, ``{Iterative methods for concave programming},'' in
  \emph{Preference, production, and capital}.\hskip 1em plus 0.5em minus
  0.4em\relax Cambridge University Press, 1989, pp. 135--148.

\bibitem{IvoBabuska73}
I.~Babu\v{s}ka, ``The finite element method with penalty,'' \emph{Mathematics
  of Computation}, vol.~27, no. 122, pp. 221--228, 1973.

\bibitem{Papadopoulos98}
P.~Papadopoulos and J.~M. Solberg, ``{Recent Advances in Contact Mechanics A
  Lagrange multiplier method for the finite element solution of frictionless
  contact problems},'' \emph{Mathematical and Computer Modelling}, vol.~28,
  no.~4, pp. 373--384, 1998.

\bibitem{seiler2012enriching}
M.~Seiler \emph{et~al.}, ``Enriching coarse interactive elastic objects with
  high-resolution data-driven deformations,'' in \emph{Proceedings of the ACM
  SIGGRAPH/Eurographics Symposium on Computer Animation}, 2012, pp. 9--17.

\bibitem{seiler2014data}
------, ``Data-driven simulation of detailed surface deformations for surgery
  training simulators,'' \emph{IEEE Transactions on Visualization and Computer
  Graphics}, vol.~20, no.~10, pp. 1379--1391, Oct 2014.

\bibitem{seiler2014efficient}
M.~U. Seiler \emph{et~al.}, ``{Efficient Transfer of Contact-Point Local
  Deformations for Data-Driven Simulations},'' in \emph{Workshop on Virtual
  Reality Interaction and Physical Simulation}.\hskip 1em plus 0.5em minus
  0.4em\relax The Eurographics Association, 2014.

\bibitem{rappel2016bayesian}
\BIBentryALTinterwordspacing
H.~Rappel \emph{et~al.}, ``Bayesian inference for the stochastic identification
  of elastoplastic material parameters: Introduction, misconceptions and
  additional insight,'' \emph{CoRR}, vol. abs/1606.02422, 2016. [Online].
  Available: \url{http://arxiv.org/abs/1606.02422}
\BIBentrySTDinterwordspacing

\bibitem{Hauseux2017}
\BIBentryALTinterwordspacing
P.~Hauseux \emph{et~al.}, ``Accelerating monte carlo estimation with
  derivatives of high-level finite element models,'' \emph{Computer Methods in
  Applied Mechanics and Engineering}, pp.~--, 2017. [Online]. Available:
  \url{http://www.sciencedirect.com/science/article/pii/S0045782516313470}
\BIBentrySTDinterwordspacing

\end{thebibliography}
\end{document}